\documentclass{IEEEtran}
\IEEEoverridecommandlockouts
\usepackage{cite}
\usepackage{amsmath,amssymb,amsfonts}
\usepackage{algorithmic}
\usepackage{graphicx}
\usepackage{textcomp}
\usepackage{xcolor}
\def\BibTeX{{\rm B\kern-.05em{\sc i\kern-.025em b}\kern-.08em
    T\kern-.1667em\lower.7ex\hbox{E}\kern-.125emX}}
\usepackage{url}

\usepackage{tikz}
\usepackage{xcolor}
\usetikzlibrary{patterns}
\usepackage{algorithm}
\usepackage{algorithmic}
\usetikzlibrary{arrows.meta,positioning} 
\usepackage{thmtools}
\usepackage{thm-restate}
\usepackage{amssymb}
\usepackage{microtype}
\usetikzlibrary{decorations.pathreplacing}

\usepackage{pgfplots}
\pgfplotsset{compat=1.12}
\usepgfplotslibrary{patchplots}
\usepgfplotslibrary{statistics}
\usepgfplotslibrary{external} 
\usepgfplotslibrary{fillbetween}
\usepackage{subfig}
\usetikzlibrary{external}
\usepackage{pgfplotstable}

\pgfplotsset{select coords between index/.style 2 args={
    x filter/.code={
        \ifnum\coordindex<#1\fi
        \ifnum\coordindex>#2\fi
    }
}}

\definecolor{colorpyC0}{RGB}{31, 119, 180} 
\definecolor{colorpyC1}{RGB}{255, 127, 14} 
\definecolor{colorpyC2}{RGB}{44, 160, 44} 
\definecolor{colorpyC3}{RGB}{214, 39, 40} 
\definecolor{colorpyC4}{RGB}{148, 103, 189} 
\definecolor{colorpyC5}{RGB}{140, 86, 75} 
\definecolor{colorpyC6}{RGB}{227, 119, 194} 
\definecolor{colorpyC7}{RGB}{127, 127, 127} 
\definecolor{colorpyC8}{RGB}{188, 189, 34} 
\definecolor{colorpyC9}{RGB}{23, 190, 207} 

\newcommand{\qedsymbol}{\hfill $\blacksquare$}
\newcommand{\pf}{\textsc{Pf}} 
\newcommand{\focs}{\textsc{Focs}}
\newcommand{\opt}{\textsc{opt}}
\newcommand{\mpc}{\textsc{Mpc}}
\newcommand{\yds}{\textsc{Yds}}
\newcommand{\avr}{\textsc{Avr}}
\newcommand{\soa}{\pf_{\textsc{oa}}}
\newcommand{\sopt}{\pf_{\opt}}
\newcommand{\woa}{w_{\textsc{oa}}}
\newcommand{\wopt}{w_{\opt}}
\newcommand{\oa}{\textsc{oa}}
\newcommand{\edf}{\textsc{Edf}}

\newtheorem{thm}{Theorem}
\newtheorem{lem}[thm]{Lemma}
\newtheorem{cor}{Corollary}
\newtheorem{defn}{Definition}

\begin{document}

\title{Relating Electric Vehicle Charging to Speed Scaling with Job-Specific Speed Limits\\
\thanks{This research is conducted within the SmoothEMS met GridShield project subsidized by the Dutch ministries of EZK and BZK (MOOI32005).}
}

\author{\IEEEauthorblockN{Leoni Winschermann, Marco E. T. Gerards, Antonios Antoniadis, Gerwin Hoogsteen, Johann Hurink}
\IEEEauthorblockA{\textit{Dept. of Electrical Engineering, Mathematics and Computer Science, University of Twente, Enschede, the Netherlands} \\
\{l.winschermann, m.e.t.gerards, a.antoniadis, g.hoogsteen, j.l.hurink\}@utwente.nl\vspace{-10pt}}
}

\author{
\IEEEauthorblockN{Leoni Winschermann, Antonios Antoniadis, Marco E. T. Gerards, Gerwin Hoogsteen, Johann Hurink}\\
\IEEEauthorblockA{Department of Electrical Engineering, Mathematics and Computer Science\\
University of Twente\\
Enschede, The Netherlands\\
\{l.winschermann, a.antoniadis, m.e.t.gerards, g.hoogsteen, j.l.hurink\}@utwente.nl}
}

\maketitle

\begin{abstract}
Due to the ongoing electrification of transport in combination with limited power grid capacities, efficient ways to schedule the charging of electric vehicles (EVs) are needed for the operation of, for example, large parking lots. Common approaches such as model predictive control repeatedly solve a corresponding offline problem. 
    
In this work, we first present and analyze the Flow-based Offline Charging Scheduler (\focs), an offline algorithm to derive an optimal EV charging schedule for a fleet of EVs that minimizes an increasing, convex and differentiable function of the corresponding aggregated power profile. To this end, we relate EV charging to processor speed scaling models with job-specific speed limits. We prove our algorithm to be optimal and derive necessary and sufficient conditions for any EV charging profile to be optimal.
    
Furthermore, we discuss two online algorithms and their competitive ratios for a specific class objective functions. In particular, we show that if those algorithms are applied and adapted to the presented EV scheduling problem, the competitive ratios for Average Rate and Optimal Available match those of the classical speed scaling problem. Finally, we present numerical results using real-world EV charging data to put the theoretical competitive ratios into a practical perspective.
\end{abstract}


\section{Introduction}
    Due to the on-going electrification of transport in combination with limited power grid capacities~\cite{2014EisingSmartGridsRisksIntegrationEnergyandTransportation} and synchronization effects~\cite{2010TuritsynRobustBroadcastCommControlofEVCharging}, 
    efficient ways to schedule the charging of electric vehicles (EVs) are needed for the operation of, for example, large parking lots. In practice, however, individual vehicles come with uncertainty in their availability and energy demand~\cite{2019VecchioMayTheForceMoveYouRolesActorsofInformationSharingDevices}. To bridge this information gap, model predictive control (\mpc) can be applied~\cite{2021VanKriekingePeakShavingCostMinMPCUniandBiDirectionalEVs}. Such \mpc\ frameworks introduce a (predictive) model to the scheduler that based on all information available at the current moment in time derives a control action for the next time step. Basic examples for such models are predictions based on historical data (e.g.,~\cite{2024IreshikaUncertaintiesInMPC4DecentrallizedAutonomousDSMofEVs}), or the introduction of deterministic charging guarantees of the form that everyone receives $x$ units of energy within $y$ hours (e.g.,~\cite{2023WinschermannPHCguaranteesvetobuttons}). 
    Another possible model is centered around the prediction of fill-levels that dictate the targeted aggregated power profile~\cite{2018SchootUiterkampFillLevelPredictioninOnlineValleyFillingAlgsforEVcharging}.
    The resulting planning may either be updated periodically, for example every 15 minutes, or rescheduling may occur based on events, for example the arrival or (early) departure of an EV. 
    As a result, \mpc s repeatedly solve an offline problem. This problem is characterized by the EVs' arrival times, departure times, energy demand and EV-specific maximum charging rates. EVs can charge simultaneously and charging of a single EV may be preempted. 
    
    To account for the limited grid capacity and a quadratic relation between charging powers and energy losses, one natural objective in EV scheduling problems is to minimize the sum of squares of the aggregated power profile of for example a parking lot hosting multiple EVs.
    EV scheduling problems with this objective 
    naturally reduce to (processor) speed scaling problems with job-specific speed limits. Hereby, as opposed to the classical model, multiple jobs may run simultaneously. In speed scaling, tasks are scheduled on a processor within their respective availability such that a (typically increasing and strictly convex) function of the processing speed is minimized. One such function may correspond to the $\ell_2$-norm, a well-studied objective function in both processor scheduling and energy research. 
    Speed scaling problems without speed limits are well-studied, with the YDS algorithm being one of the core approaches \cite{1995YaoSchedulingModelforReducedCPUEnergyYDSalg}. Already before YDS, \cite{1981VizingRussianYDS} studied the same problem and came up with a similar solution as early as 1981
    . An extension of YDS considering continuous speed limits for the aggregated speed profile is given by \cite{2017AntoniadisContSpeedScalingWithVariability}. Another variant considers job-specific speed functions and uses a maximum flow formulation to find an optimal solution for both a single-processor and multi-processor setup \cite{2017ShiouraMachineSpeedScalingbyAdaptingMethodsforConvexOptimizationwithSubmodularConstraints}. \cite{2011ZhangOptSpeedScalingAlgsSpeedChangeConstraints} investigate a model where changes in global speed are associated with additional cost.
    However, to the best of our knowledge, the use case with job-specific speed limits, as is relevant to EV scheduling with EV-specific maximum charging powers, has not yet been studied. Here, job-specific speed limits correspond to the maximum charging powers of the individual EVs. 
    As summarized in recent reviews by for example~\cite{2024ElghanamOptimizationEVSchedulingReview,2024AlAlwaschOptScheduleEVReview, 2024SinghEVchargingSchedulingReview}, centralized EV scheduling problems like the one discussed here are typically solved to optimality by formulating them as mathematical programs (e.g., linear or quadratic programs) or approximated with heuristics. The overview by~\cite{2024ElghanamOptimizationEVSchedulingReview} illustrates that the applied methods to solve mathematical programs and the used heuristics are usually generic methods that do not consider the explicit problem structure, or constructively derive an optimal schedule.
    
    In this work, we make the following contributions to both the understanding of the offline and online versions of the EV scheduling problem sketched above. 
    \begin{itemize}
        \item We discuss the relation between the classical speed scaling model and the extension based on EV scheduling.
        \item We present and analyze a novel offline algorithm to derive an optimal charging schedule for a fleet of EVs, minimizing the integral of an increasing and strictly convex function of the aggregated speed profile.
        \item We derive necessary and sufficient conditions for solutions to the offline problem to be optimal.
        \item We derive competitive ratios for two natural algorithms for the online problem known from speed scaling without job-specific speed limits. In particular, we show average rate scheduling to be directly applicable to the extended model and to have a competitive ratio of $2^{\alpha-1}\alpha^{\alpha}$ for a given $\alpha$-dependent objective function. The second algorithm considered is an adapted version of Optimal Available as previously studied for speed scaling without job-specific speed limits. We show it to be $\alpha^{\alpha}$ competitive for the same objective function. Both competitive ratios match those for the corresponding algorithms applied to speed scaling without job-specific speed limits.
        \item We put those theoretical results in perspective using real-world data to empirically quantify the performance of average rate and Optimal Available for an EV scheduling case.
        \item We prove that given the aggregated speed profile of a feasible solution, there exists no deterministic online scheduling rule that reliably finds a feasible solution following the given profile.
    \end{itemize}
    
    The remainder of the paper is organized as follows. Section~\ref{sec:problemStatement} formally describes the extended speed scaling model. After that, in Section~\ref{sec:offlineDSL} we analyze the offline problem and present the \textit{Flow-based Offline Charging Scheduler} (\focs), which is an offline algorithm that uses maximum flows to compute an optimal solution for the given problem. Then, Section~\ref{sec:onlineDSL} extends the problem to online scheduling, considering the competitive ratios of two natural algorithms. Furthermore, we prove that given the aggregated speed profile of a feasible solution, there exists no deterministic online scheduling rule that reliably finds a feasible solution that follows the given profile. Finally, we compare theoretical competitive ratios to empirical results in Section~\ref{sec:numericalEmpricialCompetitiveRatiosDSL} using numerical experiments based on real-world EV charging data. Section~\ref{sec:conclusionDSL} presents the conclusion of the paper.

\section{Problem statement}
\label{sec:problemStatement}

In this section, we describe the considered speed scaling models for processor scheduling with and without job-specific speed limits. Note that the used notation follows scheduling convention to emphasize the relation to classical results. In sentences that discuss energy, we caution the reader to carefully consider the context since the term refers to two related but different concepts in respectively processor scheduling and EV research.

The \emph{Deadline-Based Speed-Scaling with Speed Limits (DSL)} problem is defined as follows. Consider a set $\mathcal{J}:=\{1,\dots n\}$ of jobs that have to be scheduled on a speed-scalable processor. Each job $j\in \mathcal{J}$ is characterized by its \emph{workload} $p_j\in\mathbb{R}_{\ge 0}$,  \emph{release time} $r_j\in \mathbb{R}_{\ge 0}$,  \emph{deadline} $d_j\ge r_j$ as well as a \emph{job-specific} \emph{speed limit} $\ell_j\in \mathbb{R}_{\ge 0}$. A \emph{schedule} is given by a function $s: \mathbb{R}_{\ge 0}\rightarrow (\mathbb{R}_{\ge 0})^n$, such that $s(t)$ is a vector describing at what speed each job is processed at time $t$. Let $s_j(t)$ be the $j^{th}$ entry of that vector. 
\begin{defn}\label{defn:DSLfeasibleSched}
A schedule for a given set of jobs $\mathcal{J}$ is said to be \emph{feasible for DSL} if 
\begin{enumerate}
    \item[(i)] every job $j\in \mathcal{J}$ is fully processed within $[r_j,d_j)$, i.e., $\int_{r_j}^{d_j}s_j(t)\,dt \ge p_j$,
    \item[(ii)] $s_j(t) = 0$ for $t\notin [r_j,d_j)$, and 
    \item[(iii)] each job respects its speed limit, i.e., $s_j(t)\le \ell_j$ for all $t\in\mathbb{R}_{\ge 0}$.
\end{enumerate}  
\end{defn}
Note that jobs may be preempted and (in contrast to the classical speed scaling model) run simultaneously.
\begin{defn}\label{defn:pf}
For a schedule $s$ let $\pf_s:\mathbb{R}_{\ge 0}\rightarrow \mathbb{R}_{\ge 0}$ be the \emph{speed profile} of $s$ defined by $\pf_s(t) = \sum_{j}s_j(t)$.
\end{defn}
For any schedule $s$, we consider an objective function
\begin{align} 
	F(s) = \int_{0}^{\infty } \bar F (\pf _{s}(t)) \,dt \label{eq:generalizedFunctionObjective}
\end{align}
of the aggregated speed profile to quantify the \emph{intensity} of schedule $s$. The aim of the optimization problem discussed in this paper is to minimize this function $F$. Here, function $\bar F$ as used in the definition of $F$ in \eqref{eq:generalizedFunctionObjective} is a strictly convex and increasing function. 
Note that in power applications, energy losses are quadratically correlated with the speed profile, and for dynamic voltage and frequency scaling that relation is cubic. Therefore, a natural choice for objective function $F$ is the \emph{energy consumption} of a schedule, given by
\begin{align}
	E(s) = \int_{0}^\infty \left(\pf_s(t)\right)^\alpha \,dt, \label{eq:energyFunctionObjective}
\end{align}
where $\alpha>1$ is a constant.  

In the following we assume that all considered DSL instances are \emph{feasible}, i.e., they satisfy
\begin{align}\label{eq:workFeasiblyFitsAvailAssumption}
    p_j \leq \ell_j \left( d_j - r_j\right) \hspace{10pt} \forall j \in \mathcal{J}.
\end{align}

DSL is closely related to the \emph{Deadline-Based Speed-Scaling (DS)} problem. The only difference is that inputs to DS omit the speed limits, and at most one job can run at any given time. More formally, the following conditions apply to a feasible DS schedule.
\begin{defn} \label{defn:DSfeasibleSched}
A schedule is said to be \emph{feasible for DS} if 
\begin{enumerate}
    \item[(i)] every job $j\in \mathcal{J}$ is fully processed within $[r_j,d_j)$, i.e., $\int_{r_j}^{d_j}s_j(t)\,dt \ge p_j$, 
    \item[(ii)] $s_j(t) = 0 $ for $t\notin [r_j,d_j)$, and
    \item[(iii)] at most one job runs at any time, i.e., $|\{j | s_j(t)>0\}|\le 1$ for all $t\in\mathbb{R}_{\ge 0}$. 
\end{enumerate}
\end{defn}
The rest of the problem definition carries over from that of DSL. 

\section{Offline scheduling}
\label{sec:offlineDSL}
As mentioned above for the \mpc\ context, in the operational reality of EV charging DSL may be solved repeatedly. Therefore, in this section, we analyze the offline DSL problem. In particular, we analyze the relation between feasible schedules for DSL and DS (Section~\ref{sec:preliminariesofflinealgs}), derive necessary and sufficient optimality conditions for offline algorithms (Section~\ref{sec:kkt}), and introduce and analyze the Flow-based Offline Charging Scheduler (\focs), an offline algorithm that solves DSL to optimality (Sections~\ref{sec:offlineAlg} and \ref{sec:algAnalysis}).

\subsection{Preliminaries offline algorithms} \label{sec:preliminariesofflinealgs}
Let $T$ be the set of all time points that are either a release time or a deadline of a given problem instance. Formally, $T:= \{t \mid \exists j\in \mathcal{J}: t = r_j \text{ or } t = d_j\}$. We refer to the elements of $T$ as breakpoints. Let $t_1<t_2<\dots <t_{m+1}$ be the sorted elements of $T$, whereby $1\leq m\leq 2n-1$, and intervals $\{[t_i,t_{i+1})\mid i = 1,\dots m\}$ partition the time between the earliest release time and latest deadline into subintervals. We refer to those subintervals as \emph{atomic intervals}, and denote atomic interval $[t_i, t_{i+1})$ as $M_i$. Hereby, the length $t_{i+1}-t_i$ of interval $M_i$ is denoted as $|M_i|$.
Furthermore, we introduce an operator $L$ for the combined length of a set of atomic intervals, i.e., if $\mathcal{T}$ is a set of indices, then $L(\mathcal{T}) = \sum_{i\in \mathcal{T}}|M_i|$. We denote the set of all indices corresponding to atomic intervals by $\mathcal{M} = \{1,\dots,m\}$.

Given that the close relationship between DSL and DS plays a central role in our results it is useful to derive notation for the inputs, outputs and algorithms for each problem: 
\begin{defn}
Let  $\mathcal{I}_\text{DSL}$ (resp. $\mathcal{I}_\text{DS})$ be the set of all possible inputs to DSL (resp. DS), with $I\in \mathcal{I}_\text{DSL}$ being described as $ I = \langle \vec{r}, \vec{d}, \vec{p}, \vec{\ell}\rangle$ (resp. $I\in \mathcal{I}_\text{DS}$, with $ I = \langle \vec{r}, \vec{d}, \vec{p} \rangle$) where $\vec{r},\vec{d},\vec{p} \text{ and } \vec{\ell}$ refer to a vector of the respective release times, deadlines, processing volume and (if applicable) speed limits describing the job set $\mathcal{J}$. We say that an instance $I = \langle \vec{r},\vec{d},\vec{p},\vec{\ell}\rangle \in \mathcal{I}_\text{DSL}$ \emph{augments} an instance $I' = \langle \vec{r'},\vec{d'},\vec{p'}\rangle \in\mathcal{I}_{DS}$ if and only if $\vec{r}=\vec{r'}$, $\vec{d}=\vec{d'}$ and $\vec{p}=\vec{p'}$. We express this as $a(I) = I'$ and call $I$ and $I'$ \emph{corresponding}. Note that function $a:\mathcal{I}_\text{DSL}\rightarrow\mathcal{I}_\text{DS}$ is \emph{not} one-to-one.
\end{defn}

We note that a feasible schedule for an input $I'\in\mathcal{I}_\text{DS}$ is also feasible for input $I=a(I')\in\mathcal{I}_\text{DSL}$ if and only if it satisfies all job-specific speed limits.

\begin{lem}\label{lem:DSL2DS}
Any feasible schedule $s$ for an instance $I=\langle \vec{r},\vec{d},\vec{p},\vec{\ell} \rangle \in\mathcal{I}_\text{DSL}$ of  the DSL problem, can be transformed into a feasible schedule $s'$ for the corresponding augmented instance $I' = a(I) \in\mathcal{I}_\text{DS}$, such that both schedules have the same speed profile (i.e., $\pf_s = \pf_{s'}$ and therefore also $F(s)=F(s')$
).\
\end{lem}
\textit{Proof of Lemma~\ref{lem:DSL2DS}:}
Consider a schedule $s$ that is feasible for instance $I=\langle \vec{r},\vec{d},\vec{p},\vec{\ell}\rangle \in\mathcal{I}_\text{DSL}$, with associated job set $\mathcal{J}$. We show how to transform $s$ into a schedule $s'$ such that $\pf_s=\pf_{s'}$ and $s'$ is feasible for the corresponding augmented instance $I'=a(I)$. To this end, we consider atomic intervals $M_i=[t_i,t_{i+1})\in \mathcal{M}$ separately.

Schedule $s'$ is obtained by simply scheduling within each interval $M_i$ and for each job $j\in\mathcal{J}$ an amount of processing volume equal to $\int_{t_i}^{t_{i+1}}s_j(t) dt$. All these volumes are scheduled within $M_i$ according to Earliest Deadline First (\edf) and with the same speed profile $\pf_s$.

By the definition of $s'$ it is straightforward that $\pf_s=\pf_{s'}$ and therefore it remains to argue that $s'$ is a feasible schedule for $I'$. Indeed, properties $(i)$ and $(ii)$ as defined in Definition~\ref{defn:DSfeasibleSched} hold because for any atomic interval $M_i$ it is the case that $\int_{t_i}^{t_{i+1}}s_j(t) \,dt = \int_{t_i}^{t_{i+1}}s'_j(t) \, dt$ and furthermore by definition the interior of $M_i$ contains no release time or deadline. Property $(iii)$ follows directly by the definition of \edf.
\qedsymbol

\begin{figure}[]
    \centering
    \begin{tabular}{c|c c c c}
         $j$ & $r_j$ & $d_j$ & $p_j$ & $\ell_j$ \\  \hline
         1 & 0 & 1 & 1 & 2 \\
         2 & 0 & 2 & 2 & 2 
    \end{tabular}
    \\ \vspace{6 pt}
    \includegraphics{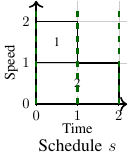}
    \includegraphics{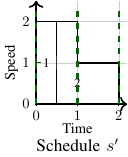}
    \caption{Illustration of transformation applied in proof of Lemma~\ref{lem:DSL2DS}.}
    \label{fig:DSL2DS}
\end{figure}
Figure~\ref{fig:DSL2DS} illustrates the method applied in the proof for instances $I=\langle(0,0),(1,2),(1,2),(2,2)\rangle$ and $I'=a(I)=\langle(0,0),(1,2),(1,2)\rangle$. Here, breakpoints $t_0$, $t_1$, and $t_2$ are indicated with green dashed lines. 

\begin{figure}[]
    \centering
    \begin{tabular}{c|c c c c}
         $j$ & $r_j$ & $d_j$ & $p_j$ & $\ell_j$ \\  \hline
         1 & 0 & 2 & 2 & 1 \\
         2 & 1 & 2 & 2 & 2 
    \end{tabular}

    \includegraphics{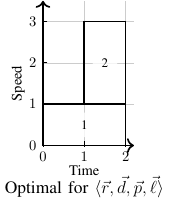}    
    \includegraphics{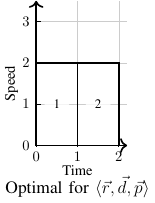}
    \caption{Example of an instance where the optimal speed profile for DSL instance $I$ differs from the optimal speed profile of the augmented DS instance $a(I)$ under objective function \eqref{eq:energyFunctionObjective} with $\alpha = 2$.}
    \label{fig:YDSinfeasibleDSLnotation}
\end{figure}

\begin{figure}[]
    \centering
    \begin{tabular}{c|c c c c}
         $j$ & $r_j$ & $d_j$ & $p_j$ & $\ell_j$ \\  \hline
         1 & 0 & 1 & 1 & 2 \\
         2 & 0 & 2 & 2 & 2 
    \end{tabular}

    \includegraphics{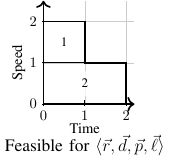}
    \includegraphics{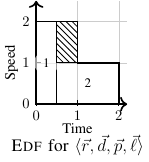}
    \caption{Example of an instance where \edf\ results in an infeasible schedule for a DSL instance even when respecting speed limits.}
    \label{fig:EDFinfeasible}
\end{figure}

Figure~\ref{fig:YDSinfeasibleDSLnotation} gives an example showing that the converse statement to that of Lemma~\ref{lem:DSL2DS} is not true. It shows an instance $I = \langle \vec{r},\vec{d},\vec{p},\vec{\ell} \rangle$, the corresponding instance $a(I)$, and a feasible schedule $s'$ for DS instance $a(I)$, for which no feasible schedule for $I$ with the exact same speed profile \pf$_{s'}$ exists. Moreover, even if given a speed profile \pf$_s$ corresponding to a feasible schedule $s$, \edf\ does not necessarily result in a feasible schedule, even if it respects job-specific maximum speeds. See Figure~\ref{fig:EDFinfeasible} for an example of this phenomenon. 

The above lemma  implies the following corollary.
\begin{cor} \label{cor:optDSscheduleSmallerEqOptDSLsched}
	The intensity of an optimal schedule $s$ for DSL instance $I = \langle \vec{r},\vec{d},\vec{p},\vec{\ell} \rangle$ is at least as high as the intensity of an optimal schedule $s'$ for the corresponding DS instance $I' = \langle \vec{r},\vec{d},\vec{p}\rangle$, i.e., $F(s) \geq F(s')$.
\end{cor}

\subsection{Optimality conditions}\label{sec:kkt}
In this section, we give a convex programming formulation for the considered offline EV scheduling problem to derive necessary and sufficient optimality conditions. To this end, we extend the mathematical program given in \cite{2007BansalSpeedScalingManageEnergyAndTemperature}. 

First, we introduce some additional notation for the offline model. 
For a given atomic interval $M_i$, we denote by $J(i)$ the jobs that are available in interval $M_i$, i.e., 
\[
    J(i) =\{ j \in \mathcal{J} | (r_j \leq t_i) \land (t_{i+1} \leq d_j)\}.
\]
Similarly, for a job $j\in \mathcal{J}$ we define by $J^{-1}(j)$ the set of indices $i$ for which job $j$ is available in interval $M_i$.
Finally, as atomic intervals are in general not unit-sized, we introduce maximum work limits $p_{i,j}^{\max} = \ell_j |M_i|$ per job $j$ and interval $M_i$.

As decision variables, let $p_{i,j}$ be the work scheduled for job $j$ during atomic interval $M_i$, i.e., $p_{i,j} = \int_{t_i}^{t_{i+1}}s_j(t)\,dt$ where schedule $s$ is yet to be determined. Given those $p_{i,j}$ values, a schedule $s$ follows naturally by scheduling job $j$ at speed $\frac{p_{i,j}}{|M_i|}$ throughout $M_i$. Therefore, decision variables $p_{i,j}$ naturally correspond to a discretized EV charging schedule where all $s_j(t)$ are \emph{step functions}, i.e., the speed between two breakpoints is constant. 
Moreover, note that by Jensen's inequality~\cite{1906JensenInequality} 
\begin{align}
    \bar F\left(\frac{\int_{t_i}^{t_{i+1}}\pf_s(t)\,dt}{|M_i|} \right) |M_i|  &\leq \int_{t_i}^{t_{i+1}}\bar F(\pf_s(t))\,dt \label{eq:jensen}
\end{align}
for any schedule $s$ and atomic interval $M_i$. Furthermore, by additivity of the integral we may consider the intensity function $F$ per atomic interval. Therefore, the right-hand-side of~\eqref{eq:jensen} is $M_i$'s contribution to the intensity of the schedule according to schedule $s$. Furthermore, the left-hand-side of~\eqref{eq:jensen} gives the intensity of a schedule that has constant speed $\frac{\int_{t_i}^{t_{i+1}}\pf_s(t)\,dt}{|M_i|}$ throughout $M_i$. 
As by strict convexity of $\bar F$ the optimal speed profile of an instance is unique, this implies that the optimal speed profile is constant within atomic intervals.
As we do not consider vehicle-to-grid applications in this work, we require that $p_{i,j}\geq 0$. Together with the feasibility conditions for DSL introduced in Section~\ref{sec:problemStatement}, we summarize the mathematical model for DSL as follows:
\begin{subequations}\label{eq:MIP}
\begin{align}
    \sum_{i\in J^{-1}(j)} p_{i,j} &\geq p_{j}  &\forall j\in\mathcal{J} \label{eq:MIPnoENS}\\
    p_{i,j} &\geq 0& \forall j \in \mathcal{J}, i\in J^{-1}(j) \hspace{3pt}\label{eq:MIPnonnegLoad}\\
    p_{i,j} &\leq p_{i,j}^{\max} & \forall j \in \mathcal{J}, i\in J^{-1}(j). \label{eq:MIPspeedlimit}
\end{align}
\end{subequations}
Note that these constraints are the same as those used by \cite{2007BansalSpeedScalingManageEnergyAndTemperature} (up to notation), extended by Inequality~(\ref{eq:MIPspeedlimit}), which models the job-specific speed limits.

From a grid perspective, the aggregated power level resulting from an EV schedule is of interest. For a given schedule, the average aggregated speed in atomic interval $M_i$ is given by $\frac{\sum_{j\in J(i)} p_{i,j}}{|M_i|}$.

Next, we consider the KKT conditions corresponding to the problem. Generally, for a convex program 
\begin{align}
\min \text{ } & \phi(x) &\nonumber\\
\text{s.t. } & \displaystyle \psi_{k}(x) \leq 0 & k = 1,\dots ,N \nonumber
\end{align}
with differentiable functions $\psi_k$ are expressed using the KKT multipliers $\lambda_k$ associated with $\psi_k$. These necessary and sufficient conditions for optimality of solutions $x$ and $\lambda$ \cite{2004BoydVandenbergheConvexOptimization} are 
\begin{subequations}\label{eq:KKTgeneral}
\begin{align}    
    \psi_k(x) &\leq 0 & k = 1,\dots , N \\
    \lambda_k &\geq 0 & k = 1,\dots , N \label{eq:KKTlambda}\\
    \lambda_k \psi_k(x) &= 0 & k = 1,\dots , N \label{eq:KKTcomplementarySlackness}\\
    \nabla \phi(x) + \sum_{k=1}^N\lambda_k \nabla \psi(x) &= 0. \label{eq:KKTgradientCondition}
\end{align}
\end{subequations}
In this section, we consider the general form $F$ of the objective function as defined in \eqref{eq:generalizedFunctionObjective}.
Applying \eqref{eq:KKTgeneral} to the discretized formulation of DSL (\ref{eq:MIP}), and introducing dual variables denoted by $\delta_j$ for (\ref{eq:MIPnoENS}), 
$\gamma_{i,j}$ for (\ref{eq:MIPnonnegLoad}), and $\zeta_{i,j}$ for (\ref{eq:MIPspeedlimit}), KKT condition \eqref{eq:KKTgradientCondition} leads to
\begin{align}
	0 &= \nabla F\left(s\right)\nonumber \\ 
	&+ \sum_{j=1}^n \delta_j \nabla \left(p_j - \sum_{i \in J^{-1}(j)} p_{i,j}\right)\label{eq:gradientNabla} \\
	&- \sum_{i \in \mathcal{M}} \sum_{j \in J(i)} \gamma_{i,j} \nabla p_{i,j}\nonumber \\
	&- \sum_{i \in \mathcal{M}} \sum_{j \in J(i)} \zeta_{i,j} \nabla \left(p_{i,j}^{\max} - p_{i,j}\right).\nonumber
\end{align}
First, note that $F(s)$ is the integral of a strictly convex and increasing function $\bar F$ of the aggregated speed profile and that the optimal speed profile is constant within atomic intervals and assumes value $\frac{\sum_{j\in \mathcal{J}p_{i,j}}}{|M_i|}$. On top of that, by additivity of the integral in $F$, decision variable $p_{i,j}$ only effects the part of the schedule's intensity associated with atomic interval $M_i$. In the following, we use this fact to gain insights on speed levels of atomic intervals based on an analysis of the components of the gradient in~\eqref{eq:gradientNabla} with respect to $p_{i,j}$.
Note that the component of this gradient that corresponds to the partial derivative with respect 
to $p_{i,j}$ is
\begin{align}\label{eq:nablacomponent}
	0 = \frac{\partial F
	}{\partial p_{i,j}} -\delta_j -\gamma_{i,j} + \zeta_{i,j}.
\end{align}

We analyze condition \eqref{eq:nablacomponent} for components corresponding to partial derivatives with respect to $p_{i,j}$, where job 
$j\in J(i)$. We consider three cases in our analysis. 

First, consider $0 < p_{i,j} < p_{i,j}^{max}$. In this case job $j$ charges in interval $i$, but not at full power. 
Complementary slackness (see (\ref{eq:KKTcomplementarySlackness})), now implies that $p_{i,j}\gamma_{i,j} = 0$ and $(p_{i,j} - p_{i,j}^{\max}) \zeta_{i,j} = 0$. 
In the considered case, this implies that $\gamma_{i,j} = \zeta_{i,j} = 0$. Therefore, (\ref{eq:nablacomponent}) simplifies to 
\begin{align}
	& 0 = -\delta_j + \frac{\partial F}{\partial p_{i,j}} \nonumber\\
	\iff & \delta_j = \frac{\partial F}{\partial p_{i,j}}.
	\label{eq:KKTanalysisDelta}
\end{align}
This shows that the dual variable $\delta_j$ is the derivative of the intensity function $F$ with respect to $p_{i,j}$. 
Since $\delta_j$ does not depend on $i$, if there is another atomic interval $i'\in J^{-1}(j)$ where $0<p_{i',j}<p_{i',j}^{\max}$, we have $\frac{\partial F}{\partial p_{i,j}} = \frac{\partial F}{\partial p_{i',j}}$.
Substituting with the definition of $F$ (cf.~\eqref{eq:generalizedFunctionObjective}) and using additivity of the integral and the fact that the aggregated speed outside of $M_i$ is independent of $p_{i,j}$, this yields
\begin{align*}
	\frac{\partial \left(\int_{t_i}^{t_{i+1}} \bar F(\pf_s(t))\,dt\right)}{\partial p_{i,j}} &=  \frac{\partial \left(\int_{t_{i'}}^{t_{i'+1}} \bar F(\pf_s(t))\,dt\right)}{\partial p_{i',j}}.
\end{align*}
Note that previously we derived that the aggregated speed within atomic intervals is constant. Therefore, we may substitute $\pf_s(t)$ by $\frac{\sum_k p{i,k}}{|M_i|}$ on domain $M_i$ of the integral, and find that 
\begin{align*}
	|M_i|\frac{\partial\left(\bar F\left(\frac{\sum_{k}p_{i,k}}{|M_i|}\right)\right)}{\partial p_{i,j}} &= |M_{i'}|\frac{\partial\left(\bar F\left(\frac{\sum_{k}p_{i',k}}{|M_{i'}|}\right)\right)}{\partial p_{i,j}}. 
\end{align*}
Here, we consider the numerator to be a composition of functions. By applying the chain rule, we find that 
\begin{align*}
	\frac{\partial\left(\bar F\right)}{\partial p_{i,j}}\left(\frac{\sum_{k}p_{i,k}}{|M_i|}\right) &= \frac{\partial\left(\bar F\right)}{\partial p_{i,j}}\left(\frac{\sum_{k}p_{i',k}}{|M_{i'}|}\right),
\end{align*}
where we now evaluate the partial derivatives of $\bar F$ in $\frac{\sum_{k}p_{i,k}}{|M_{i}|}$ and $\frac{\sum_{k}p_{i',k}}{|M_{i'}|}$ respectively.
By strict convexity of $\bar F$, we conclude that $\frac{\sum_{k}p_{i,k}}{|M_{i}|}=\frac{\sum_{k}p_{i',k}}{|M_{i'}|}$, i.e., the aggregated speeds of any such two intervals $M_i$ and $M_{i'}$ where job $j$ charges at a rate strictly between 0 and its power limit is the same. 

Next, consider the case where $0 = p_{i,j} < p_{i,j}^{\max}$. Complementary slackness gives $\zeta_{i,j} = 0$, leaving us with 
\begin{align}
	& 0 = -\delta_j + \frac{\partial F}{\partial p_{i,j}} - \gamma_{i,j} \nonumber\\
	\iff & \gamma_{i,j} = -\delta_j + \frac{\partial F}{\partial p_{i,j}}.
	\label{eq:KKTanalysisGamma}
\end{align}
Using non-negativity of $\gamma_{i,j}$ (see (\ref{eq:KKTlambda})), it follows that $\frac{\partial F}{\partial p_{i,j}} \geq \delta_j$. As above, $\delta_j$ is independent of $i$ and characterizes $\frac{\partial F}{\partial p_{i',j}}$ for intervals with index $i'$ where $0<p_{i',j}<p_{i,j}^{\max}$. 
Using the same arguments as above and the fact that the derivative of a strictly convex and increasing function is increasing, 
we conclude that the aggregated speed during interval $M_i$ where by assumption job $j$ does not process at positive speed is at least as high as during intervals where job $j$ does process at a (positive) speed below its maximum. 

Lastly, consider the case where $0 < p_{i,j} = p_{i,j}^{\max}$. Complementary slackness gives us $\gamma_{i,j} = 0$, leaving us with 
\begin{align}
	& 0 = -\delta_j + \frac{\partial F}{\partial p_{i,j}} + \zeta_{i,j} \nonumber\\
	\iff & \zeta_{i,j} = \delta_j - \frac{\partial F}{\partial p_{i,j}}.
	\label{eq:KKTanalysisZeta}
\end{align}
Applying similar reasoning as in the previous cases and considering that the signs in the right hand sides of (\ref{eq:KKTanalysisGamma}) and (\ref{eq:KKTanalysisZeta}) are reversed, we conclude that the speed in any interval $M_i$ where job $j$ is executed at maximum speed, is at most as high as during intervals where $j$ is available and is either processed at a (positive) power below its maximum, or is available and not processed at all.

From the above analysis, the following necessary and sufficient conditions for a schedule to be optimal follow:
\begin{enumerate}
    \item[\textsc{Kkt1}] The aggregated speed in all intervals where $j$ is scheduled but does not reach its speed limit is the same.
    \item[\textsc{Kkt2}] The aggregated speed in intervals where $j$ could, but does not run is at least as high as in intervals where $j$ actually runs.
    \item[\textsc{Kkt3}] The aggregated speed in intervals where $j$ runs at maximum speed is smaller or equal than in intervals where $j$ runs below its speed limit. 
\end{enumerate}
The first two conditions are similar to those derived by Bansal, Kimbrel and Pruhs, whereas the last results from the addition of job-specific speed limits. 

In Section~\ref{sec:optimality} we show that the output of the \focs\ algorithm introduced in Section~\ref{sec:alg} is a feasible schedule that satisfies said conditions. For such a schedule, we can solve the system (\ref{eq:KKTanalysisDelta}), (\ref{eq:KKTanalysisGamma}) and (\ref{eq:KKTanalysisZeta}), proving optimality of the derived primal solution. 

\subsection{Offline algorithm using flows} \label{sec:offlineAlg}
In this section, we present an iterative offline algorithm to determine an optimal schedule for DSL instances, minimizing the integral of an increasing and strictly convex function of the aggregated speed profile. First, note that due to the convexity of the objective function and the finite number of release times and deadlines, the aggregated speed profile of any optimal solution is a step function. Moreover, the aggregated speed within any atomic interval $M_i$ is constant for such a schedule. Similarly to YDS, the algorithm presented here uses the notion of critical intervals. These intervals are exactly those intervals that in an optimal solution require the highest aggregated power. Formally, these intervals are defined as follows.
\begin{defn}[Critical intervals] \label{def:crit}
    An atomic interval $M_i$ is \textit{critical} if for any optimal schedule $s$ its average aggregated speed $\frac{1}{|M_i|}\int_{t_i}^{t_{i+1}} \pf_{s}(t)\,dt$ is larger or equal to the average aggregated speed $\frac{1}{|M_{i'}|}\int_{t_{i'}}^{t_{i'+1}} \pf_{s}(t)\,dt$ for any $i'\in \mathcal{M}$.
\end{defn}

Note that there may be multiple critical (atomic) intervals. Furthermore, one major difference with critical intervals as defined for YDS is that jobs do not have to be fully contained within a (set of) critical interval(s) in order to be scheduled there. This difference with YDS follows from the job-specific speed limits. The speed profile that YDS assigns to what they call a critical interval when solving DS instances is not necessarily feasible in DSL, the setting with speed limits (see e.g., Figure~\ref{fig:YDSinfeasibleDSLnotation}). 
Compared to YDS, determining critical intervals and their power level is more involved. In the algorithm presented in Section~\ref{sec:alg}, determining critical intervals is based on the computation of multiple maximum flows. 
To be able to compute the flows and to keep track of the developments over the iterations of the proposed algorithm, we follow the DSL 
notation introduced so far, and introduce some additional notation.

\subsubsection{Flow formulation.}\label{sec:flow}
\begin{figure}[]
    \centering
    \includegraphics{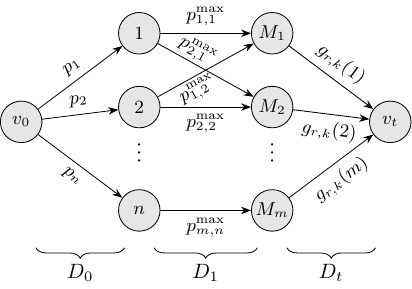}
    \caption{Schematic of flow network structure of DSL.}
    \label{fig:flowSchematic}
\end{figure}

For the proposed algorithm, we use a network $G = (V,D)$. The network is initialized as follows (see also Figure~\ref{fig:flowSchematic}).
The vertex set V consists of source and sink vertices $v_0$ and $v_t$, as well as two sets of vertices representing job vertices and atomic interval vertices respectively, i.e., $V = \{v_0,v_t\}\cup\mathcal{J}\cup\{M_i|i\in \mathcal{M}\}$.
Furthermore, the edge set $D$ consists of the union of the following three sets:
\begin{align*}
    D_0 &= \{(v_0,j) | j\in \mathcal{J}\} \\
    D_1 &= \{(j,M_i) | j\in \mathcal{J}, i\in J^{-1}(j)\} \\
    D_t &= \{(M_i,v_t) | i\in \mathcal{M}\}
\end{align*}
with respective edge capacities
\begin{align*}
    c_{u,v} = \begin{cases}
        p_v & \text{if } u = v_0, \ v \in \mathcal{J} \\
        p_{i,u}^{max} & \text{if } u \in \mathcal{J}, \ v = M_i,\ i\in J^{-1}(u)\\
        g_{r,k}(u) & \text{if } u \in \mathcal{M}, \ v = v_t
    \end{cases}.
\end{align*}
Note that the function $g_{r,k}$ is not defined yet. The algorithm works with rounds (indexed by $r$), each of which executes iterations (indexed by $k$). Intuitively, $g_{r,k}$ is a lower bound on the flatness of the aggregated speed profile. It varies over the execution of the algorithm, and is discussed in more detail in Section~\ref{sec:alg}. 

Given a flow $f$ in network $G$, we denote the flow value as $|f|$ and call an edge $(u,v)$ \textit{saturated} if $f(u,v) = c_{u,v}$.
Note that a flow in $G$ corresponds to an EV schedule. Here, a job $j$ is scheduled to process $f(j,M_i)$ units of work in interval $M_i$, or equivalently an EV $j$ charges $f(j,M_i)$ in interval $M_i$
and the capacities on edges in $D_1$ model the job-specific speed limits. 
Furthermore, any flow for which the edges in $D_0$ are saturated corresponds to a feasible EV charging schedule 
and the flow through $D_t$ models the aggregated speed in the atomic intervals of the charging schedule corresponding to $f$. Note that the capacities and flow through $D_t$ are expressed in terms of the aggregated work processed. 
By normalizing for the length of each atomic interval, we can deduce the aggregated speed profile.
Based on this correspondence, we may use the network structure to not only derive a feasible, but an optimal schedule for objective function $F(\sum_{j\in \mathcal{J}}s_j(t))$.

\subsubsection{Algorithm formulation.} \label{sec:alg}
In the following, we use network $G$ defined in Section~\ref{sec:flow} to derive an iterative algorithm that gives an optimal schedule and power profile for (aggregated) EV charging with objective function $F(s) = \int_0^{\infty} \bar F(\pf_s(t))\,dt$ where $\bar F$ is strictly convex and increasing (cf.~\eqref{eq:generalizedFunctionObjective}).

Before going into detail, we provide some intuition and a rough overview of the workings of the algorithm. 
Intuitively, the edge set $D_0$ can be interpreted as the processing work of the jobs. For any feasible DSL schedule, those demands have to be met. The flow through edge set $D_1$, on the other hand, is what we are trying to determine: the schedule itself. For any interval node $M_i$, the incoming flow corresponds to the load scheduled in that interval. In particular, $f(j,M_i)$ is the work processed for job $j$ in interval $M_i$. 
Whereas the capacities of edges in $D_0$ and $D_1$ are determined by the instance, edge capacities of edges in $D_t$ are not. However, the flow through $D_t$ directly corresponds to the value of the objective function. Therefore, the algorithm presented in this section defines edge capacities for $D_t$ such that they are a lower bound on the highest aggregated speed contributing to the objective function, i.e., a lower bound on the outgoing flow of nodes $M_i$ where $M_i$ is a critical interval. If given those capacities, we find a maximum flow that saturates all edges in $D_0$, we have found a feasible solution with this maximum speed, and use this to determine the partial schedule for any critical interval $M_i$. This partial schedule corresponds to the incoming flow at each such node $M_i$. Else, we adapt the lower bound and repeat the process until we do find such a maximum flow and partial schedule. 

\begin{figure}[]
    \centering
    \includegraphics{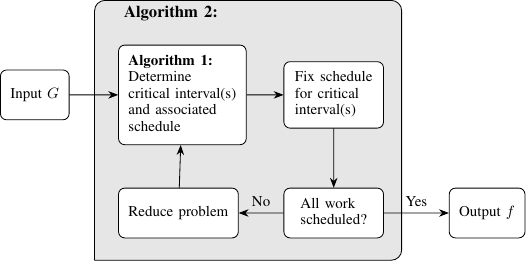}
    \caption{Schematic overview of \focs.}
    \label{fig:algSchematic}
\end{figure}

In Figure~\ref{fig:algSchematic} we provide a rough outline of an algorithm that exploits the bottleneck function of the critical intervals.
In both the algorithm formulation and analysis, we distinguish between iterations and rounds of the algorithm. In Figure~\ref{fig:algSchematic}, a new round starts every time that Algorithm~\ref{alg:round} is called. To determine a (set of) critical interval(s) (see Definition~\ref{def:crit}), we may require multiple iterations in which we adapt the lower limit. Given the dynamic nature of this lower limit, we denote it as $g_{r,k}$ where $r$ and $k$ denote the current round and iteration respectively. 
At the end of a round, we have determined a (set of) critical interval(s). We determine the schedule for those intervals to be the incoming flow at the corresponding interval nodes. For non-critical intervals, there is no schedule yet. Their schedules will be determined in the next rounds. In that fashion, we will construct a schedule for the entire instance. 
To keep track of what has yet to be scheduled, we introduce the notion of \textit{active intervals}. At the beginning of a round, interval $M_i$ is active if it has not yet been scheduled (i.e., has not yet been critical) in previous rounds. Let $\mathcal{M}_a$ be the set of indices of active intervals, which we initialize to be all atomic intervals, i.e., $\mathcal{M}_a = \mathcal{M}$. 

The first iteration of the first round goes as follows. Given that we have to schedule a certain amount of energy and that the objective function is increasing, the most optimistic lower bound on the aggregated power is a constant profile over all intervals. Therefore, we initialize the capacities of the edges in $D_t$ by
\begin{align}
    g_{1,1}(i) = \frac{\sum_{j=1}^n p_j}{L\left(\mathcal{M}_a\right)}|M_i| \ \forall i \in \mathcal{M}_a,\ \nonumber
\end{align}
which is the aggregated energy charged in atomic interval $M_i$ given that in all intervals the same aggregated charging power is used, and all energy requirements are met. In that way, the edge capacities of $D_t$ act as lower bounds to the highest aggregated power level. They are dynamic and will be increased over iterations. Given the capacities, we determine a maximum flow $f_{1,1}$ for this instance. If the flow value $|f_{1,1}|$ of $f_{1,1}$ is $\sum_{j=1}^n p_{j}$, we have found a feasible schedule, and all active intervals are critical. If not, then there is at least one non-saturated edge $(M_i,v_t)$ with $M_i$ an active interval. We call the intervals corresponding to such edges \textit{subcritical}. Note that those intervals will not be critical in this round. We therefore temporarily remove them from the set of active intervals and add them to what we call the collection of \textit{parked intervals} $\mathcal{M}_p$. At the beginning of each round, this collection is initialized to be empty. This is the end of the first iteration. 

From here, we structurally increase the edge capacities of edges in $D_t$ and again compute a maximum flow until all edges in $D_0$ are saturated, and we find a feasible EV schedule. To this end, first note that after the first iteration, 
\begin{align}
    \sum_{j=1}^n c_{v_0,j} - |f_{1,1}| = \sum_{j=1}^n p_j - |f_{1,1}| > 0 \nonumber
\end{align}
if there were subcritical intervals. In particular, this means that there are jobs $j$ for which additional work still needs to be scheduled. Among the interval-vertices, the only candidates for additional flow are those vertices $M_i$ for which edge $(M_i,v_t)$ was saturated in $f_{1,1}$, i.e., the remaining active intervals.
Keeping the objective in mind, we therefore proportionally increase the capacities at the remaining active intervals to
\begin{align}
    g_{1,2}(i) = g_{1,1}(i) + \frac{\sum_{j=1}^n p_j - |f_{1,1}|}{L\left(\mathcal{M}_a\right)}|M_i| \ \forall i \in \mathcal{M}_a. \nonumber
\end{align}

We repeat this process until we find a flow with flow value $\sum_{j=1}^n p_{j}$. Such a flow leads to a feasible EV schedule for which the maximum aggregated power is minimal. Say this happens after $K_1$ iterations. 
The remaining active intervals in that iteration make up the set of critical intervals in the corresponding round. In Figure~\ref{fig:algSchematic}, this case corresponds to the first time we leave the box of Algorithm~\ref{alg:round} and move on to fix parts of the schedule we aim to compute.

We generalize the steps discussed thus far to an arbitrary round $r$ and iteration $k$ with $1\leq k <K_r-1$ where $K_r$ is the number of iterations in round $r$. This yields: 
\begin{align*}
    g_{r,1}(i) &= \frac{\sum_{j=1}^n c_{v_0,j}}{L\left(\mathcal{M}_a\right)}|M_i| \ &\forall i \in \mathcal{M}_a \\
    g_{r,k+1}(i) &= g_{r,k}(i) + \frac{\sum_{j=1}^n c_{v_0,j} - |f_{r,k}|}{L\left(\mathcal{M}_a\right)}|M_i| \ &\forall i \in \mathcal{M}_a
\end{align*}
given that flow $f_{r,k}$ is the maximum flow in round $r$ and iteration $k$, and that between iterations active intervals and flow networks are updated. We end the round when we find a maximum flow with flow value $\sum_{j=1}^n c_{v_0,j}$.

After each round $r$, we fix the part of the schedule associated with the critical interval(s) (top right box in Figure~\ref{fig:algSchematic}) to correspond to the flow incoming at the respective (critical) interval nodes, and reduce the remainder of the problem by constructing a new network $G_{r+1}$ (bottom left box in Figure~\ref{fig:algSchematic}) as follows. First, we exploit the acyclic topology of the network to define a flow $f_r|_{M_r^*}$ of the determined maximum flow $f_r$, where $M_r^* = \{i \in \mathcal{M} | M_i \text{\ is critical in round\ }r \}$ is the set of indices of critical intervals and 
\begin{align*}
    f_r|_{M_r^*}(M_i,v_t) &= 
        \begin{cases}
            f(M_i,v_t)  & \text{if } i \in M_r^* \\
            0       & \text{otherwise}
        \end{cases} \\
    f_{r}|_{M_r^*}(j,M_i) &= 
        \begin{cases}
            f(j,M_i)  & \text{if } i \in M_r^* \\
            0       & \text{otherwise}          
        \end{cases} \\
    f_r|_{M_r^*}(v_0,j) &= \sum_{i\in J^{-1}(j)} f_r|_{M_r^*}(j,M_i).
\end{align*} 
Note that this definition backpropagates flow from the sink to the source. Intuitively, $f_r|_{M_r^*}$ denotes the flow that goes through critical intervals. In the YDS-sense, $f_r|_{M_r^*}(v_0,j)$ is the critical load of job $j$ in round $r$. 
Now, $G_{r+1}$ is the network obtained by removing edges $(M_i,v_t)$ with $i\in I^*_r$ from $G_r$, and updating edge capacities to be $c_{u,v} - f_r|_{M_r^*}(u,v)$. From here, we start the next round of the algorithm and initialize a new flow $f_{r+1,1}$. 
For convergence, between iterations within a round we similarly construct $G_{r,k+1}$ based on the subcritical flow $f_{r,k|\mathcal{M}_p}$. Alternatively, we can require for $k \geq 1$ that we initialize flow $f_{r,k+1}$ with $f_{r,k}$ and augment it to a maximum flow using for example shortest augmenting path algorithms. 

The optimal flow output by the algorithm is $f = \sum_r f_r|_{M_r^*}$. Implicitly, we use that augmenting paths in future rounds will not reshuffle the already determined subschedule induced by the critical intervals. We will come back to that in Lemma~\ref{lemma:Isolation}. For more information about augmenting paths, and their relation to maximum flows, please refer to e.g., \cite{1972EdmondsKarpMaxFlowAlgs}. 

\begin{figure}[]
    \centering
    \begin{tabular}{c|c c c c}
         $j$ & $r_j$ & $d_j$ & $p_j$ & $\ell_j$ \\  \hline
         1 & 0 & 3 & 2 & 2 \\
         2 & 1 & 2 & 2 & 2 \\
    \end{tabular}
    \\ \vspace{6 pt}
    \includegraphics{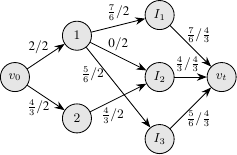}
    \includegraphics{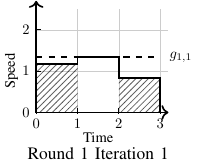}
    \\ \vspace{6pt}
    \includegraphics{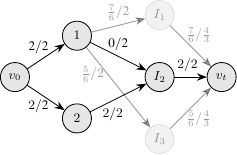}
    \includegraphics{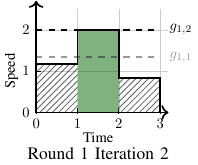}
    \\ \vspace{6pt}
    \includegraphics{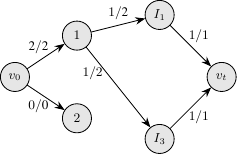}
    \includegraphics{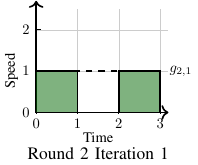}
    \\ \vspace{6 pt}
    \includegraphics{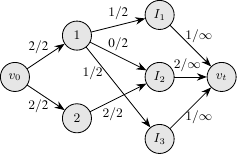}
    \includegraphics{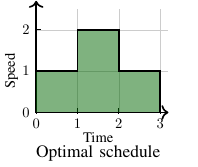}
    \caption{Intermediate states of \focs \ for an example instance, tracked over rounds and iterations.}
    \label{fig:flowExWithProfile}
\end{figure}

The algorithm to derive a feasible schedule within a round is summarized in Algorithm~\ref{alg:round}. This algorithm is then embedded in the global algorithm (Algorithm~\ref{alg:alg}) described in this section, outputting a flow $f$ corresponding to an optimal EV charging schedule. We refer to this algorithm as \textit{Flow-based Offline Charging Scheduler} (\focs). To illustrate, Figure~\ref{fig:flowExWithProfile} displays both the flow and aggregated power profile of an example instance over the rounds and iterations of the algorithm. Here, the first three flows display $f_{r,k}$, whereas the last flow is the optimal flow $f$. In the respective power profiles corresponding to the flow-induced schedules, shaded intervals are parked, and solid green intervals are critical. In general, maximum flows are not unique. To illustrate that, the first flow is deliberately chosen such that the flow through $(1,I_1)$ differs from that through $(1,I_3)$. Note how the optimal power profile in the bottom graph is a sum of the green components at the end of each round of the algorithm. Note that Step~\ref{alg:alg:repeat} in Algorithm~\ref{alg:alg} can be reformulated as a recursion by calling \focs($G_r$).

\begin{algorithm}[] 
\begin{algorithmic}[1]
 \renewcommand{\algorithmicrequire}{\textbf{Input:}}
 \renewcommand{\algorithmicensure}{\textbf{Output:}}
    \caption{\textsc{Round}}\label{alg:round}
    \REQUIRE $
                        G_r, \
                        r, \
                        \mathcal{M}_a
    $\\
    \ENSURE $
                        \text{feasible flow } f_r, \
                        \text{critical sets } M_r^*
    $\\
    \STATE{Initialize: }$
                        \mathcal{M}_p = \varnothing, \ 
                        k = 0, \
                        G_{r,k} = G_r
    $
        \STATE 
            \label{alg:rep:defineG}
            $c_{M_i,v_t} = g_{r,k}(i) \ \forall i \in \mathcal{M}_a$
        \STATE Determine a maximum flow $f_{r,k}$
        \STATE 
            $\mathcal{M}_p = \mathcal{M}_p \cup \{i\in \mathcal{M}_a | i \text{ subcritical in } f_{r,k}\}$ 
        \STATE $\mathcal{M}_a = \mathcal{M}_a \setminus \mathcal{M}_p$
        \IF{$|f_{r,k}| = \sum_{j=1}^n c_{v_0,j}$} \label{alg:round:ifcondition}
            \RETURN $f_r = f_{r,k},\ M_r^* = \mathcal{M}_a$
        \ELSE
            \STATE $G_{r,k+1} = G_{r,k}$ with capacities reduced by $f_{r,k|\mathcal{M}_p}$ and vertices $M_i$ removed for subcritical $M_i$
            \STATE $k = k + 1$ and repeat from Step~\ref{alg:rep:defineG}        
        \ENDIF
\end{algorithmic}
\end{algorithm}

\begin{algorithm}[] 
\begin{algorithmic}[1]
 \renewcommand{\algorithmicrequire}{\textbf{Input:}}
 \renewcommand{\algorithmicensure}{\textbf{Output:}}
    \caption{\textsc{Focs}}\label{alg:alg}
    \REQUIRE $
                        G
    $
    \ENSURE $
                        \text{optimal flow } f
    $
    \STATE {Initialize: }$
                        \mathcal{M}_a = \mathcal{I}, \ 
                        \mathcal{M}_p = \varnothing, \ 
                        r = 0, \
                        G_r = G , \ 
                        f
    $
        \STATE \label{alg:alg:callRep}
            $f_r$, $M_r^* = $ \textsc{Round}($G_r,r,\mathcal{M}_a$) 
                \STATE $\mathcal{M}_a = \mathcal{M}_a \setminus M_r^* $
                \STATE $f = f + f_r|_{M_r^*}$
                \STATE $G_{r+1} = G_r$ with capacities reduced by $f_r|_{M_r^*}$ and vertices $M_i$ removed for $i\in M_r^*$ 
                \STATE $r = r + 1$
        \IF{$\mathcal{M}_a = \emptyset$} \label{alg:alg:ifcondition} 
            \RETURN $f$
        \ELSE 
                \STATE Repeat from Step~\ref{alg:alg:callRep} \label{alg:alg:repeat}
        \ENDIF
\end{algorithmic}
\end{algorithm}

\subsection{Algorithm analysis}\label{sec:algAnalysis}
In this section, we analyze \focs, the algorithm presented above. In particular, Section~\ref{sec:properties} shortly discusses its time complexity and properties, after which its optimality is proved in Section~\ref{sec:optimality}.

\subsubsection{Properties and time complexity.}\label{sec:properties}
In this section, we discuss some properties and lemmas that apply to the flow model and algorithm. In particular, we establish some building blocks that enable us to prove optimality of the algorithm in Section \ref{sec:optimality}. 

\begin{lem} \label{lem:focsfeasibleterminates}
    If an instance has a feasible schedule, \focs\ terminates and outputs a feasible schedule. Its time complexity is bound by $\mathcal{O}(n^2\mu)$ where $\mathcal{O}(\mu)$ is the time complexity of the used maximum flow algorithm.
\end{lem}
\textit{Proof of Lemma~\ref{lem:focsfeasibleterminates}: }    
    First, we argue that \focs\ terminates and analyze its time complexity. We do this by arguing that the number of iterations in Algorithm~\ref{alg:round}, and the number of rounds in Algorithm~\ref{alg:alg} calling Algorithm~\ref{alg:round} are finite. 
    
        Any feasible schedule $s$ can be directly translated to a feasible flow for $G_r$ by sending $\int_{t_i}^{t_{i+1}}s_j(t)\,dt$ units of flow through edge $(j,M_i)$ for each job $j$ and $i\in J^{-1}(j)$. The flows through edges in $D_0$ and $D_t$ follow directly by flow conservation. 
        Therefore, there exists a maximum flow for the input to Algorithm~\ref{alg:round} which saturates all edges in $D_0$.
        The algorithm (and therefore the current round) finishes once the if condition in Step~\ref{alg:round:ifcondition} is satisfied, i.e., if we find a maximum flow that saturates all edges outgoing of sink node $v_0$ when using edge capacities $g_{r,k}$ for the network. If Step~\ref{alg:round:ifcondition} is false, there exists at least one subcritical interval $M_i$ for which the flow through $(M_i,t)$ is strictly below capacity $g_{r,k}(i)$. In each iteration, at least one such interval is removed from the network, until $g_{r,k}$ is increased sufficiently to find a maximum flow satisfying the if condition. The number of intervals is finite, and therefore the number of iterations in Algorithm~\ref{alg:round} is finite. In particular, we can bound this number to at most $2n - 1$ iterations per round, based on the fact that the number $m$ of atomic intervals is bound by the number of jobs $j$, implying that $m \leq 2n$. Note that there are efficient algorithms available to solve maximum flow problems, e.g., \cite{1956FordFulkersonMaxFlows,1972EdmondsKarpMaxFlowAlgs,1970DinitzMaxFlows,1974KarzanovMaxFlowsPreflows}. Furthermore, a comprehensive overview of traditional polynomial time maximum flow algorithms is given by \cite{1988GoldbergMaxFlows}. Denoting their time complexity by $\mu$, we find that Algorithm~\ref{alg:round} has a time complexity of $\mathcal{O}(n \mu)$.  
        
        As at the end of each round at least one interval is critical and therefore removed from $\mathcal{M}_a$, the finite number of intervals implies that the if condition in Step~\ref{alg:alg:ifcondition} of Algorithm~\ref{alg:alg} is satisfied after at most $2n-1$ rounds, and hence \focs\ terminates. This implies that the time complexity of \focs\ is bound by $\mathcal{O}(n^2\mu)$.

    Note that the time complexity for the EV charging setting may be reduced further by exploiting the underlying structure of EV charging schedules, and by considering the decrease in network size over the rounds of the algorithm. In particular, we may initialize the flow of any iteration with the flow found in the previous iteration of the same round. 
    Furthermore, there are maximum flow algorithms that are cubic in the number of nodes \cite{1988GoldbergMaxFlows}. As the largest flow network that is considered in \focs\ (the network in the initial round) has $n+m+2 \leq 3n +2$ nodes, a rough upper bound of the time complexity of maximum flows in \focs\ is given by $\mu \leq n^3$.
    
    Finally, feasibility of the output follows from the defined edge capacities of the considered network. Following the order of DSL feasibility conditions listed in Definition~\ref{defn:DSLfeasibleSched}, we conclude that:
    \begin{itemize}
        \item[(i)] Every job $j$ is fully processed within its availability since by flow conservation the exact amount of work done per job within its availability is the flow through edge $(v_0,j)$. The algorithm only terminates once that edge is saturated, i.e., if $f(v_0,j) = p_j$. 
        \item[(ii)] An edge $(j,M_i)$ is in $D_1$ if and only if $j\in J(i)$. Therefore, assuming the default value is zero, all decision variables $e_{i,j}$ for which $i\in \mathcal{M}\setminus J^{-1}(j)$ are zero, implying that the speed of $j$ outside its availability is zero.
        \item[(iii)] Each job respects its speed limit by the capacities defined for edges in $D_1$.
    \end{itemize}
    Therefore, the output of \focs\ is feasible.
\qedsymbol

Next, we extend on the concept of work-transferability as described by \cite{2017AntoniadisContSpeedScalingWithVariability} to integrate job-specific speed limits.
\begin{defn}[Work-transferability]
    If for a given schedule and atomic intervals $M_i$ and $M_{i'}$ there exists a job $j\in J(i) \cap J(i')$ such that $p_{i,j} > 0$ and $p_{i',j} < p_{i',j}^{\max}$, we state that the \textit{work-transferable} relation $i\rightarrow i'$ holds. Furthermore, let $\twoheadrightarrow$ be the transitive closure of $\rightarrow$.
\end{defn}

Intuitively, if we have work-transferability from one atomic interval $i$ to another atomic interval $i'$, then we can transfer some work that was scheduled during $i$ to $i'$. In EV charging terms this implies that we can advance or delay some charging from one period in time to another. Applying the concept to flows, we can make the following statement.

\begin{lem}[Work-transferability in flows] \label{lemma:WT}
    For a given schedule and atomic intervals $M_i$ and $M_{i'}$, we have $i\rightarrow i'$ if and only if there exists a path $(M_i,j,M_{i'})$ in the residual graph corresponding to the schedule, where $j\in \mathcal{J}$. Similarly, we have $i\twoheadrightarrow i''$ if and only if in the residual network corresponding to the schedule there exists an ($M_i$-$M_{i''}$)-path through interval and job vertices only. 
\end{lem}
\textit{Proof of Lemma~\ref{lemma:WT}:}
We show only the first statement as the extension follows naturally using concatenations of paths. Assume that $i\rightarrow i'$. Then there exists a job $j$ such that $j\in J(i) \cap J(i')$ with $p_{i,j} > 0$ and $p_{i',j} < p_{i',j}^{\max}$. The former implies that edge $(M_i,j)$ exists in the residual graph. As $c_{j,M_{i'}} = p_{i',j} < p_{i',j}^{\max}$, edge $(j,M_{i'})$ is in the residual graph. This proves existence of path $(M_i,j,M_{i'})$ in the residual graph. 

For the opposite direction, assume the existence of a path $(M_{i},j,M_{i'})$. Since $j\in \mathcal{J}$, we know the edge capacity $c_{j,M_{i'}}$ in the original network to be $p_{i',j}^{\max}$. 
The existence of the edge in the residual graph implies that for the flow going through this edge which is defined by the schedule to be $p_{i',j}$, we have $p_{i',j} < p_{i,j}^{\max}$.
Furthermore, existence of edge $(M_i,j)$ in the residual graph indicates positive flow through $(j,M_i)$ in the original network
, implying $p_{i,j} > 0$. From the presence of both edges, it follows that $j\in J(i) \cap J(i')$, proving that $i\rightarrow i'$.
\qedsymbol

\begin{figure}[]
    \centering
    \includegraphics{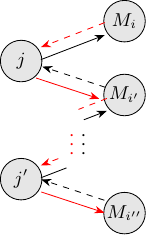}
    \caption{Work-transferability represented in flows.}
    \label{fig:wtFlow}
\end{figure}

Figure~\ref{fig:wtFlow} illustrates the concept of work-transferability. Here, dashed edges are those that are not in the original network, but might be present in the residual network. Lemma~\ref{lemma:WT} translates work-transferability to the existence of (in the figure) red paths in the residual network. 

Next, we consider two lemmas that have a more direct relation to the algorithm. 

\begin{figure}[]
    \centering
    \includegraphics{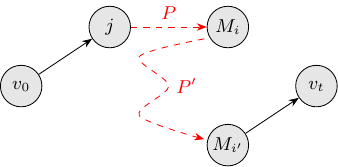}
    \caption{Illustration of augmenting path in proof of Lemma~\ref{lemma:Isolation}.}
    \label{fig:isolationProof}
\end{figure}

\begin{lem}[Isolation of critical intervals] \label{lemma:Isolation}
    If $M_i$ is a critical interval in round $r$ and if the round consists of multiple iterations whereby $M_{i'}$ was subcritical in one of those iterations, then there is no work-transferable relation between $i$ and $i'$ in the schedule corresponding to the flow at the end of round $r$, i.e., $i\not \twoheadrightarrow i'$ with respect to flow $f_r$.
\end{lem}
\textit{Proof of Lemma~\ref{lemma:Isolation}:}
We prove the lemma by constructing an augmenting path (see Figure~\ref{fig:isolationProof}). Assume in round $r$ interval $M_{i'}$ was parked in iteration $k$ and let $f_{r,k}$ be the flow at the end of iteration $k$. 
Since $M_{i'}$ is subcritical, we have $|f_{r,k}| < \sum_{j=1}^n p_j$, implying that for the next iteration the lower bound $g_{r,k}$ will be increased to
\begin{multline*}
        g_{r,k+1}(i'') = g_{r,k}(i'') +  \frac{\sum_{j=1}^n c_{v_0,j} - |f_{r,k}|}{L\left(\mathcal{M}_a\right)}|M_{i''}| \\ \forall i''  \in \mathcal{M}_a.
\end{multline*}
By criticality of $M_i$, the interval is active at the end of the iteration, implying $g_{r,k+1}(i) > g_{r,k}(i)$. Furthermore, criticality implies that there is no iteration in this round where $M_i$ is subcritical. Combing those facts, the flow through $(M_i,t)$ increases in iteration $k+1$ compared to iteration $k$.
This is only possible if there is a job $j$ such that $(v_0,j)$ is not saturated and there exists a ($j$-$M_i$)-path $P$ in the residual graph.
Furthermore, note that since $M_{i'}$ is being parked in iteration $k$, edge $(M_{i'},t)$ is not saturated and therefore exists in the residual graph.
Now, assume $i\twoheadrightarrow i'$. By Lemma~\ref{lemma:WT}, there exists an ($M_i$-$M_{i'}$)-path $P'$ that passes only through job and interval vertices. This implies that $P'' = (v_0,P,P',v_t)$ exists in the residual graph and contains an ($v_0$-$v_t$)-path, proving existence of an augmenting path in $f_{r,k}$. This contradicts maximality of the flow, implying $i \not \twoheadrightarrow i'$.
\qedsymbol

Intuitively, this lemma says that we cannot push any charging from (high power) critical intervals to (low power) subcritical intervals. This is in line with the notion of critical intervals as introduced for the YDS algorithm, and will be a key element in the optimality proof in Section \ref{sec:optimality}. Furthermore, this particular lemma justifies that we fix the schedule of critical intervals at the end of each round. 

For the next lemma, we first introduce the notion of ranks. 

\begin{defn}[Rank]
    The rank $r(i)$ of an atomic interval $M_i$ is defined as the round $r(i)$ in which $M_i$ was critical, i.e., $i\in M_{r(i)}^*$.
\end{defn}

\begin{lem}[Monotonicity] \label{lemma:Monotonicity}
    For the schedule corresponding to output flow $f$ of algorithm \focs\ and atomic intervals $M_i$ and $M_{i'}$ where $r(i) < r(i')$, the aggregated power in $M_i$ is strictly larger than in $M_{i'}$, i.e.,
    \begin{align}
        \frac{f(M_i,t)}{|M_i|} > \frac{f(M_{i'},t)}{|M_{i'}|}\nonumber
    \end{align}
\end{lem}
\textit{Proof of Lemma~\ref{lemma:Monotonicity}:}
We prove the lemma by contradiction, whereby we consider flows at the end of rounds. Let interval $M_i$ be the lowest ranked interval such that its aggregated power level in the flow outputted by the algorithm is larger than that in intervals with rank $r(i) - 1$. As the algorithm does not change schedules at critical intervals, this already occurs at round $r(i)$ itself. Given that $M_i$ was subcritical in the previous round, the speed level for $M_i$ increased. In particular, there is a job $j$ for which the speed during interval $M_i$ increased compared to the previous round. Furthermore, we know that the flow through $(j,M_i)$ is positive in round $r(i)$, implying that edge $(M_i,j)$ is in the residual graph. However, applying Lemma~\ref{lemma:Isolation} to the previous round, the amount of work scheduled for job $j$ remains the same. Therefore, there is an interval $M_{i'}$ for which the flow through $(j,M_{i'})$ decreased compared to the previous round. As a consequence, the flow in round $r(i)$ does not saturate the edge, implying that edge $(j,M_{i'})$ is in the residual graph. Combining these findings, the path $(M_i,j,M_{i'})$ is in the residual graph, implying $i\rightarrow i'$, and thus contradicting Lemma~\ref{lemma:Isolation}.
\qedsymbol

The lemma shows that the average aggregated speed of atomic intervals is decreasing in their rank. 
Therefore, critical intervals as determined using the method presented in this paper share the monotonicity property known for YDS for corresponding DS instances. We first find those intervals with the highest intensity, and then iteratively determine the next highest speeds.  

We also note that, similarly to YDS, the power profile outputted by \focs\ is unique if the objective function is strictly convex. However, this does not necessarily apply to the schedule.

\subsubsection{Optimality proof.}\label{sec:optimality}
In this section, we prove that Algorithm~\ref{alg:alg} as described in Section~\ref{sec:alg} computes a solution $s$ that is optimal under objective function $F(s)$. To this end, we first prove some auxiliary lemmas that show compliance with the sufficient conditions derived in Section~\ref{sec:kkt}. 

\begin{lem} \label{lemma:kkt1}
    The output of Algorithm~\ref{alg:alg} complies with \textsc{Kkt1}.
\end{lem}
\textit{Proof of Lemma~\ref{lemma:kkt1}:}
    If in the final output of the algorithm there are two distinct atomic intervals $M_i$ and $M_{i'}$ such that for job $j$ we have $0 < \frac{p_{i,j}}{|M_{i}|} = \frac{p_{i',j}}{|M_{i'}|} < \ell_j $, then by definition of worktransferability  we have $i\rightarrow i'$ and $i'\rightarrow i$. By Lemma~\ref{lemma:Isolation} and the strict monotonicity in Lemma~\ref{lemma:Monotonicity}, this implies that the aggregated speed in both intervals is the same.
\qedsymbol

\begin{lem} \label{lemma:kkt2}
    The output of Algorithm~\ref{alg:alg} complies with \textsc{Kkt2}.
\end{lem}
\textit{Proof of Lemma~\ref{lemma:kkt2}:}
    Let $i\in J^{-1}(j)$ be such that $p_{i,j} = 0$ in the output of the algorithm. Assume that there is an interval $M_{i'}$ with $i\neq i'$ and $i'\in J^{-1}(j)$ for which the aggregated power in $M_{i'}$ is strictly greater than in $M_{i}$, i.e., $\frac{\sum_{j=1}^n p_{i,j}}{|M_{i}|} < \frac{\sum_{j=1}^n p_{i',j}}{|M_{i'}|}$. By Lemma~\ref{lemma:Monotonicity} we have $r(i') < r(i)$, implying by Lemma~\ref{lemma:Isolation} that $i' \not \twoheadrightarrow i$. Applying the definition of work-transferability, it follows that $p_{i',j} = 0$, proving compliance with \textsc{Kkt2}.
\qedsymbol

\begin{lem} \label{lemma:kkt3}
    The output of Algorithm~\ref{alg:alg} complies with \textsc{Kkt3}.
\end{lem}
\textit{Proof of Lemma~\ref{lemma:kkt3}:}
    Let job $j$ run at maximum speed in $M_i$ in the schedule found by \focs . Assume that there is an interval $M_{i'}$ with $i\neq i'$ and $i'\in J^{-1}(j)$, such that the aggregated speed in $M_i$ is strictly greater than in $M_{i'}$. By Lemma~\ref{lemma:Monotonicity}, we know that $r(i)<r(i')$. Therefore, by Lemma~\ref{lemma:Isolation}, there is no work-transferable relation between $i$ and $i'$ ($i\not \twoheadrightarrow i'$). From the definition of work-transferability it now follows directly that $p_{i',j} \geq p_{i',j}^{\max}$, proving compliance with \textsc{Kkt3}. 
\qedsymbol

Combining all discussed above, we conclude optimality of the algorithm output.

\begin{thm}[Optimality]\label{thm:optimality}
    For any feasible input instance, the schedule produced by Algorithm~\ref{alg:alg} is an optimal solution minimizing any convex, increasing and differentiable objective function of the aggregated output powers. 
\end{thm}
\textit{Proof of Theorem~\ref{thm:optimality}:}
    The proof follows directly from the \textsc{Kkt} conditions derived in Section~\ref{sec:kkt}, the inherent feasibility of the output and Lemmas~\ref{lemma:kkt1}--\ref{lemma:kkt3}.
\qedsymbol

To summarize, this section considered the offline DSL problem, applicable to EV scheduling in \mpc\ settings. In particular, we analyzed the relation between solutions of DSL instances and their corresponding DS instances. Furthermore, we derived necessary and sufficient optimality conditions for DSL schedules and presented an offline algorithm that determines an optimal schedule in $\mathcal{O}(n^2 \mu)$ time where $\mu$ is the complexity of an efficient maximum flow algorithm. Lastly, we provided proof of the optimality of the output of the algorithm.

\section{Online scheduling}
\label{sec:onlineDSL}
As discussed in the introduction, \mpc\ is a much-deployed framework for coordinated EV charging, especially due to its usability to bridge data gaps. However, it is an interesting and important question to ask how close to an optimal solution such frameworks can get, and in particular, how close to optimal an \mpc, or intraday controller, can get, assuming perfect knowledge on an EV's characteristics upon arrival. Note that in that case, the model-component of the \mpc\ was clairvoyant.

Also from a theoretical point of view, considering DSL in an online setting is a natural next step.
Therefore, in this section, we are interested in schedules constructed \emph{online}, i.e., schedules where jobs are released one by one, and the algorithm only gets to know their characteristics at their respective release times. 

\subsection{Preliminaries online algorithms}
We define the online variant of a job scheduling problem to be such that the existence and characteristics of jobs become known at their respective release times. 
In this section, we analyze online algorithms for DSL in terms of their respective \emph{competitive ratio}. 
\begin{defn}
    Given a deterministic algorithm \textsc{Alg} that for any DSL instance $I \in \mathcal{I}_{DSL}$ determines a feasible schedule $s^{\textsc{Alg}}(I)$, and given an optimal solution $s^*(I)$, the \emph{competitive ratio} of the algorithm is defined as
    \begin{align}
    \sup_{I\in\mathcal{I}_{DSL}} \frac{E(s^{\textsc{Alg}}(I))}{E(s^{*}(I))}. \label{eq:defCompetitiveRatio}
    \end{align}
    The definition carries over to DS instances.
\end{defn}

Two classical online approaches for DS are Average Rate (\avr) and Optimal Available (\oa)~\cite{1995YaoSchedulingModelforReducedCPUEnergyYDSalg}. Given the connection between DSL and DS, we first provide a short description of those two algorithms, before relating them to DSL.

\subsubsection{\avr.}
\avr\ for DS works in two steps.
First, upon release, job $j$ is scheduled at speed $\frac{p_j}{d_j - r_j}$ throughout its availability, i.e., each job is scheduled at the constant speed corresponding to the average speed it needs to complete its work between release time and deadline. 
The resulting schedule may not be feasible for DS, since jobs may run simultaneously, but it gives a useful initial speed profile. Therefore, in a second step, \edf\ is applied using the speed profile resulting from the initial schedule.

\subsubsection{\oa.}
\oa\ reoptimizes the remaining problem instance each time a new job is released. In particular, let $s$ be an optimal schedule for jobs $\mathcal{J} = \{1,\dots n\}$ and instance  $I = \langle \vec{r},\vec{d},\vec{p} \rangle$. Let $t'$ be the first point in time where a new job $n+1$ is released. We then define the remaining instance at point $t'$ according to schedule $s$ as $I' = \langle \vec{r'},\vec{d'},\vec{p'} \rangle$ where:
\begin{align*}
    \vec{r'}_j &= t'  \hspace{82pt} \forall j\in \mathcal{J}\cup\{n+1\} \\
    \vec{d'}_j &= d_j \hspace{80pt} \forall j\in \mathcal{J}\cup\{n+1\} \\
    \vec{p'}_j &= \begin{cases}
                    p_j - \int_{r_j}^{t'} s_j(t) \,dt & \text{if } j\in \mathcal{J}\\
                    p_{n+1} & \text{if } j=n+1.
                \end{cases}
\end{align*}
For convenience, if a job has remaining workload 0 or if its deadline is at most $t'$, we remove it from the remaining problem instance. We now determine an optimal schedule $s'$ for $I'$. The updated schedule $\bar s$ for \oa\ at time $t'$ is such that
\begin{align*}
    \bar s_j(t) = \begin{cases}
                        s_j(t) & \text{if } j\in \mathcal{J} \wedge t<t' \\
                        s'_j(t) & \text{if } j\in \mathcal{J} \wedge t\geq t' \\
                    s'_{n+1}(t) & \text{if } j = n+1.
                \end{cases}
\end{align*}
Iteratively repeated over the time horizon every time a new job is released, this results in a schedule $s_{\oa}$ for \oa. 

Note that \oa\ may be applied to either DS or DSL, with the difference being the algorithm applied in the optimization subroutine. Optimization for DS may be done by applying YDS, whereas for DSL \focs\ is a suitable optimization algorithm.

\subsection{Average Rate for DSL}
In the following, we discuss the application of \avr\ to DSL instances. As remarked earlier, \edf\ does not necessarily result in a feasible schedule for DSL instances, even if it follows a profile for which there exists a feasible schedule (see Figure~\ref{fig:EDFinfeasible}). However, since feasible schedules for DSL allow more than one job to be processed at any time and since by assumption \eqref{eq:workFeasiblyFitsAvailAssumption} holds, we can adapt \avr\ to DSL instances by skipping the last step (therefore not applying \edf) to find a feasible schedule. In other words, upon release, we schedule any job $j$ at speed $\frac{p_j}{d_j - r_j}$ for the next $d_j - r_j$ units of time.  

We analyze the performance guarantee of applying \avr\ to DSL instances by relating the resulting schedules to those resulting from applying \avr\ to DS instances as follows. 
Assume we are given an instance $I = \langle \vec{r},\vec{d},\vec{p},\vec{\ell} \rangle$. Let $s^{\text{DSL},\avr}$ be the schedule for $I$ found by \avr\ without \edf\ and let $s^{\text{DSL},*}$ be an optimal schedule. Furthermore, let $s^{\text{DS},\avr}$ be the schedule for the corresponding augmented DS instance $a(I)$ found by \avr\ with \edf, and let $s^{\text{DS},*}$ be an optimal schedule for $a(I)$. Note that \pf$_{s^{\text{DS},\avr}} = $\pf$_{s^{\text{DSL},\avr}}$ and therefore their objective values are the same. Thus,
\begin{subequations}
\begin{align}
    E(s^{\text{DSL},\avr}) 
        &= E(s^{\text{DS},\avr}) \label{eq:avrCRavrEqavr}\\
        &\leq 2^{\alpha - 1}\alpha^\alpha E(s^{\text{DS},*}) \label{eq:avrCRdsCR}\\
        &\leq 2^{\alpha - 1}\alpha^\alpha E(s^{\text{DSL},*}). \label{eq:avrCRdsSmallerEqdsl}        
\end{align}
\end{subequations}
Here, (\ref{eq:avrCRavrEqavr}) follows from the fact that the \avr\ (respectively with and without \edf) schedules $s^{\text{DS},\avr}$ and $s^{\text{DSL},\avr}$ have the same speed profile, (\ref{eq:avrCRdsCR}) follows from the known tight competitive ratio for \avr\ with \edf\ for DS instances~\cite{2007BansalSpeedScalingManageEnergyAndTemperature} and Lemma~\ref{lem:DSL2DS}, and (\ref{eq:avrCRdsSmallerEqdsl}) follows from Corollary~\ref{cor:optDSscheduleSmallerEqOptDSLsched}.

Furthermore, we can conclude that the upper bound on the competitive ratio for \avr\ without \edf\ for DSL instances is the same as the upper bound for the competitive ratio for \avr\ with \edf\ for DS instances. In particular, assume that $I' = \langle \vec{r},\vec{d},\vec{p}\rangle$ was an instance for which the inequality in (\ref{eq:avrCRdsCR}) is an equality.
Based on this, we can construct an instance $I$ for DSL such that $a(I) = I'$ by taking $\ell_j = \max_t \pf_{s^{\text{DS},*}}(t)$. Then, the energy of the respective optimal schedules of $I$ and $I'$ are the same.

\subsection{Optimal Available for DSL}
In this section, we adapt the potential function approach that Bansal, Kumar and Pruhs used to analyze the competitive ratio of \oa\ for DS instances \cite{2007BansalSpeedScalingManageEnergyAndTemperature} to analyze \oa\ for DSL instances. In particular, we show the competitive ratio of $\alpha ^\alpha$ to be tight, where $\alpha$ is the same as in the energy function defined in (\ref{eq:energyFunctionObjective}). Notably, the competitive ratio for DS and DSL instances is the same. 

First, we remark that while in each re-optimization step of \oa\ there is a unique aggregated speed profile, the schedule is not necessarily unique. Therefore, we explicitly note that throughout the upcoming analysis, we consider a fixed (arbitrary) realization of \oa.

Next, we introduce some notation, before giving the potential function, and deriving the competitive ratio.
Let $\soa(t)$ be the aggregated speed at which \oa\ runs at time $t$, and similarly let $\sopt(t)$ be the aggregated speed at which \opt\ runs at time $t$, where \opt\ is an optimal algorithm leading to an optimal (offline) schedule $s_{\opt}$. Note that these are the speed profiles as realized at the end of the time horizon and that \oa\ recomputes a schedule every time a new task is released. Therefore, for current time $t_0$, we introduce schedule $s$ and corresponding speed $\pf_s(t)$, at which \oa\ runs at time $t\geq t_0$ if no new jobs are released after $t_0$. While it holds that $\pf_s(t_0) = \soa(t_0)$, this is generally not the case for $t>t_0$, since a job may be released in interval $(t_0,t)$.

At this point $t_0$, \oa\ may be interpreted to solve a DSL instance where all tasks have the same release time $t_0$. Therefore, the resulting speed profile $\pf_s(t)$ for $t > t_0$ is a non-increasing step function. 
Let $t_1$ be the first point in time after $t_0$ where a step occurs and the speed profile changes, i.e.,
\begin{align*}
    t_1 = \sup\{t\geq t_0 | \pf_s(t) = \pf_s(t_0)\}.
\end{align*}
Applying this inductively, we find breakpoints $t_i$ for $i\geq 0$, such that $[t_i, t_{i+1})$ are the inclusion maximal intervals with constant speed profiles. As  \cite{2007BansalSpeedScalingManageEnergyAndTemperature}, we call such intervals \textit{critical}, and define $M_i := [t_i, t_{i+1})$. Note that both criticality and the intervals $M_i$ are being redefined here, as compared to their use in the offline Section~\ref{sec:offlineDSL}.

We define $\woa(t,t')$ as the work done by \oa\ in interval $(t,t']$ that is already available (and unfinished) at current time $t_0$. 
Moreover, let $\frac{\woa(t,t')}{t' - t}$ be the density of interval $(t,t']$. We note that this definition of $w$ is different from the one given by \cite{2007BansalSpeedScalingManageEnergyAndTemperature}, to account for the methods needed to solve DSL.
For critical intervals $M_i$, we can relate this work load to the speed profile by
\begin{align}
    \pf_s(t) = \frac{\woa(t_i,t_{i+1})}{t_{i+1} - t_i} \nonumber
\end{align}
for $t \in M_i$. For notational purposes, we shorten $\woa(t_i,t_{i+1})$ to $\woa(i)$ in places where $t_i$ and $t_{i+1}$ are clear from the context.

\begin{figure}[]
    \centering
    \begin{tabular}{c|c c c c}
         $j$ & $r_j$ & $d_j$ & $p_j$ & $\ell_j$ \\  \hline
         1 & 0 & 3 & 5 & 2 \\
         2 & 1 & 2 & 1 & 2 
    \end{tabular}
    \\ \vspace{6 pt}
    \includegraphics{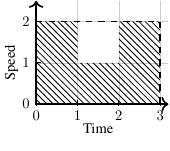}
    \caption{Example of an instance where the speed profile corresponding to available work $\wopt$ at time $t_0 = 0$ (profile of shaded area) is not a non-increasing step function. The dashed line corresponds to the optimal speed profile at the end of the time horizon.
    }
    \label{fig:WoptNotStepFunction}
\end{figure}

Similarly to $\woa$, we define $\wopt(t,t')$ to be the work done by optimal schedule \opt\ in interval $(t,t']$ that is already available at current time $t_0$. 
Note that as opposed to \oa, \opt\ is aware of tasks that will be released in the future. Therefore, the speed profile induced by $\wopt$ is not necessarily non-increasing, as shown in Figure~\ref{fig:WoptNotStepFunction}. 
The figure shows the part of the optimal speed profile corresponding to jobs that are already available. In the middle, there is a valley in the profile, where the second (not yet available) job will be scheduled.
Moreover, note that the breakpoints determined based on $\pf_s(t)$ above do not necessarily align with changes in speed in the (partial) speed profile of \opt.

Finally, we define the potential function at current time $t = t_0$ to be
\begin{align}
    \Phi (t) = \alpha \sum_{i\geq 0} \pf_s(t_i)^{\alpha - 1} \left(\woa(i) - \alpha \wopt(i)\right). \nonumber
\end{align}
Note that breakpoints $t_i$ on the right hand side of the formula depend on function input $t$. Furthermore, note that the form of the introduced potential function is similar to the potential function used by~\cite{2007BansalSpeedScalingManageEnergyAndTemperature} to derive the competitive ratio for \oa\ under DS. The main differences lie in the problem definition, and in the definition of $\woa$ and $\wopt$. The proof for DSL follows the same structure as that for DS presented by~\cite{2007BansalSpeedScalingManageEnergyAndTemperature}. However, it is not evident that the same form of potential function applies to DSL. Therefore, we work out the details in places where speed limits play a role and include the proof for DSL despite the similarities.

In the following we show that:
\begin{lem}\label{lem:potentialFunction}
    The potential function $\Phi(t)$ satisfies the following conditions:

\end{lem}
\begin{enumerate}
    \item \label{item:potentialProperty3} \textit{Boundary property:} $\Phi(t) = 0$ for any $t$ before the first release time and for any $t$ after the last deadline.
    \item \label{item:potentialProperty1} \textit{Job release and completion property:} At any time $t$ where a new task is released, or a task is completed by $s_\opt$ or $s_\oa$, the potential function is non-increasing, i.e., 
    \begin{align}\label{eq:potentialProperty1LimitFormulation}
        \lim_{t'\uparrow t} \Phi(t') \geq \lim_{t'\downarrow t} \Phi(t').
    \end{align}
    \item \label{item:potentialProperty2} \textit{General property:} For a time $t$ at which no job is released, we have
        \begin{align}
            \soa(t)^\alpha + \frac{d\Phi(t)}{dt} \leq \alpha ^\alpha \sopt(t)^\alpha. \label{eq:PFproperty2}
        \end{align}
\end{enumerate}

\cite{2010BansalSpeedScalingWeightedFlowTimePotentialFtn} have previously proven similar properties for the potential function of another speed scaling model to derive the corresponding competitive ratio. The proof of Lemma~\ref{lem:potentialFunction} uses the following lemma taken from \cite{2007BansalSpeedScalingManageEnergyAndTemperature}. We refer to the original paper for its proof.
\begin{lem}\label{lem:algebraicBansal} (\textit{Lemma 3.3 in} \cite{2007BansalSpeedScalingManageEnergyAndTemperature}\textit{): }
    Let $q, r, \delta \geq 0$ and $\alpha \geq 1$. Then $(q + \delta)^{\alpha -1}(q - \alpha r - (\alpha - 1)\delta ) - q^{\alpha -1} (q - \alpha r) \leq 0$.
\end{lem}
\textit{Proof of Lemma~\ref{lem:potentialFunction}: } We consider each property separately.

    For Property~\ref{item:potentialProperty3}, note that for any $t$ before the first release, or after the last deadline, there will be no tasks to schedule for \oa, resulting in speed $\pf_s(t_i) = 0$ for all $i\geq 0$. Hereby, $\Phi(t) = 0$ for such $t$.

    For Property~\ref{item:potentialProperty1}, we first consider the case where a new job is released. 
        Assume task $j$ is released at time $t_0$ with deadline $d_j\in M_i$, requiring an amount $x$ of work.
        This release may affect the breakpoints of the critical intervals. Note that the speed profile $\pf_s(t)$ before the release was a non-increasing step function. To add a job, we add its associated work to the latest intervals for which the job is available. We consider the addition of work in increments, i.e., we add some $x' \leq x$ of work until one of the following cases occurs:
        \begin{itemize}
            \item The speed in $M_i$ increases to the speed in $M_{i-1}$. In particular $\pf_s(t_i) = \frac{w(t_i, t_{i+1}) + x'}{t_{i+1} - t_i} = \pf_s(t_{i-1})$. 
            \item Two critical intervals $M_i$ and $M_{i+1}$ merge into one new critical interval.
            \item The speed at which job $j$ runs in interval $M_i$ reaches the job-specific speed limit $\ell_j$.
            \item The interval $M_i$ splits in two critical intervals $M_{i'}$ and $M_{i''}$. This may occur due to deadlines that do not match the already existing breakpoints, or due to speed limits being reached in parts of recently merged critical intervals.
            \item The job is completely scheduled (i.e., $x' = x$).
        \end{itemize}
        Due to the speed limits, we have to carefully keep track of the added work $x'$ so far, before adding the next part of the work, until the whole amount of work $x$ is scheduled.
        
        We start with cases that do not change the structure of critical intervals. Those are the cases where $M_i$ increases to the speed in $M_{i-1}$, where in $M_i$ job $j$ reaches its job-specific speed limit $\ell_j$, and where the remaining work can be scheduled within $M_i$ without triggering any of the other events. 
        
        By definition, \oa\ schedules all additional work $x'$ during $M_i$. Therefore, the only values associated with \oa\ that change are $\pf_s(t)$ for $t\in M_i$ and $\woa(i)$. For the optimal schedule \opt, no such claim can be made. Therefore, we denote $x_{i'}' \geq 0$ to be the work scheduled for interval $M_{i'}$ for $0 \leq i' \leq i$ where $\sum_{i'= 0}^i x_{i'}' = x'$. 
        
        We initially consider the $i^{th}$ term of the potential function separately. Speed function $\pf_s(t)$ changes by $\frac{x'}{(t_{i+1} - t_i)}$ for $t\in M_i$. To compare the change in potential function, we denote the new speed as 
        \begin{align}
            \pf_{s'}(t) = \begin{cases}
                \frac{\woa(i) + x'}{t_{i+1} - t_i} & \text{if } t\in M_i \\
                \pf_s(t) & \text{otherwise.}
            \end{cases}
        \end{align}
        We further note that $\woa(i)$ increases by $x'$, and $\wopt(i)$ increases by $x_{i}'$. That gives a total change of
        \begin{align}
            &\pf_{s'}(t_i)^{\alpha-1} \left(\woa(i) + x' - \alpha \left(\wopt(i) + x_{i}'\right)\right)\nonumber \\
            -& \pf_s(t_i)^{\alpha -1} \left(\woa(i) - \alpha \wopt(i)\right) \label{eq:PFithTerm}
        \end{align}
        for the $i^{th}$ term. In the term for $i' \in \{0,\dots,i-1\}$, the values $\woa(i')$ do not change after adding work $x'$ since the work is added to a different interval, namely $M_i$. However, $\wopt(i')$ increases by $x_{i'}'$. If we sum the change for all such $i'$, we find the following expression:
        \begin{align}
            \sum_{i'=0}^{i-1} & \left( \pf_s(t_{i'})^{\alpha - 1} \left(\woa(i') - \alpha\left(\wopt(i') + x_{i'}'\right)\right)\right) \nonumber\\ 
            &- \pf_s(t_{i'})^{\alpha - 1}\left(\woa(i') - \alpha \wopt(i')\right) \nonumber \\
            &= \sum_{i'=0}^{i-1} \pf_s(t_{i'})^{\alpha - 1} \left(-\alpha x_{i'}'\right). \label{eq:PFnonithTerm}
        \end{align}
        
        To show that (\ref{eq:potentialProperty1LimitFormulation}) holds in the release case, we denote $\Delta\Phi(t) = \lim_{t'\uparrow t} \Phi(t') - \lim_{t'\downarrow t}\Phi(t')$. We bring all terms together and conclude:
        \begin{align}
            \Delta&\Phi(t_0) \\
                &= \sum_{i'=0}^{i-1} \pf_s(t_{i'})^{\alpha -1} \left(-\alpha x_{i'}'\right) \label{eq:DeltaPFstart}\\
                    &\ + \pf_{s'}(t_i)^{\alpha - 1} \left(\woa(i) + x - \alpha \left(\wopt(i) + x_{i}'\right)\right) \nonumber\\
                    &\ - \pf_s(t_i)^{\alpha -1} \left(\woa(i) - \alpha \wopt(i)\right) \nonumber \\ 
                &\leq \pf_{s'}(t_i)^{\alpha - 1} (\woa(i) + x - \alpha (\wopt(i) + \sum_{i'=0}^{i} x_{i'}')) \label{eq:DeltaPFleqst}\\
                    &\ - \pf_s(t_i)^{\alpha -1} \left(\woa(i) - \alpha \wopt(i)\right) \nonumber \\ 
                &= \pf_{s'}(t_i)^{\alpha - 1} \left(\woa(i) + x - \alpha \left(\wopt(t_i, t_{i+1}) + x'\right)\right) \label{eq:DeltaPFreducedBansal} \\
                    &\ - \pf_s(t_i)^{\alpha -1} \left(\woa(i) - \alpha \wopt(i)\right) \nonumber \\ 
                &\leq 0. \nonumber
        \end{align}
        Here, (\ref{eq:DeltaPFstart}) follows from the derivations in (\ref{eq:PFithTerm}) and (\ref{eq:PFnonithTerm}). In (\ref{eq:DeltaPFleqst}), we made use of the fact that for all $0 \leq i' < i$, we can lower bound the speed $\pf_s(t_{i'})$ by the new speed $\pf_{s'}(t_i)$. By doing so, we can exploit that $\sum_{i'= 0}^i x_{i'}' = x'$ in the next step. Moreover, (\ref{eq:DeltaPFreducedBansal}) is exactly the case discussed in the proof of Theorem~3.4 by \cite{2007BansalSpeedScalingManageEnergyAndTemperature} where they do inequality manipulations, apply Lemma~\ref{lem:algebraicBansal} (Lemma~3.3 in \cite{2007BansalSpeedScalingManageEnergyAndTemperature}) and conclude non-positivity. For the details, we refer the reader to that work.
        
        Both for the cases where interval $M_i$ either splits in two, or merges with another, at that point the densities between old and new intervals are constant, leaving the potential function unchanged. 
        
        We conclude, that the potential function does not increase if a new task is released. 
        
        Next, consider the case where \oa\ finishes a(t least one) task at time $t$. 
        Either 
        \begin{align}\label{eq:prop1OAfinishLimitEq}
            \lim_{t'\uparrow t} \pf_s(t') = \lim_{t'\downarrow t} \pf_s(t'),
        \end{align}
        in which case the task finished strictly within the critical interval, leaving $\pf_s(t)$ unaffected and continuously reducing $\woa(t_0, t_1)$ by which (\ref{eq:potentialProperty1LimitFormulation}) holds by continuity in $t$, or the equality in (\ref{eq:prop1OAfinishLimitEq}) does not hold, in which case we transition from one critical interval to the next. Then, indices shift by one, $\pf_s(t_1)$ becomes $\pf_s(t_0)$ etc. For the potential function, the only change is the formerly first term $\alpha \pf_s(t_0)^{\alpha -1} (\woa(t_0, t_1) - \alpha \wopt(t_0, t_1))$ disappearing. However, while $\pf_s(t_0)$ remained constant, both $\woa(t_0,t_1)$ and $\wopt(t_0,t_1)$ approached zero. The first term's contribution therefore approaches zero from above as the critical interval draws to an end. The change in potential function in such a point is therefore continuous and (\ref{eq:potentialProperty1LimitFormulation}) holds.

        For the case where at time $t$, \opt\ finishes a task, we note that the potential function is independent of $\sopt(t)$, and the change in $\wopt(t_0,t_1)$ is continuous. 

    Combining the observations above, we have shown that $\Phi(t)$ satisfies Property~\ref{item:potentialProperty1}.

    Lastly, we show that $\Phi(t)$ has Property~\ref{item:potentialProperty2} by showing that 
    \begin{align}
        \soa(t)^\alpha - \alpha ^\alpha \sopt(t)^\alpha + \frac{d\Phi(t)}{dt} \leq 0. \label{eq:potentialFtnIneq}
    \end{align}
    
        We consider the working case where in the next infinitesimally small $dt$ time units no new task is released or completed by either \oa\ or \opt. 
        Furthermore, we do assume that there are tasks available, and that therefore $\soa(t_0) = \pf_s(t_0)>0$ and $\sopt(t_0)>0$ for current time $t_0$. As \oa\ runs, $\woa(t_0,t_1)$ is reduced at rate $\pf_s(t_0)$. For $i > 0$, the value of $\woa(i)$ remains unchanged. 
        For the work done by \opt, we remark once more that \opt's speed profile does not necessarily align with breakpoints $t_i$. However, it is easy to verify that if at any point $t > t_0$ the speed $\sopt(t)$ increases, at least one new task will be released at time $t$.
        Therefore, assuming that no new tasks are released in the next $dt$ units of time, we can assume $\sopt(t)$ to be a non-increasing step function over interval $[t_0,t_0 + dt)$. 
        From this, we use that $\sopt(t_0)$ is an upper bound on the rate at which \opt\ reduces $\wopt(t_0,t_1)$ throughout the next $dt$ units of time. Therefore, for (\ref{eq:potentialFtnIneq}) to hold, it suffices to show that the following inequality holds:        
        \begin{multline}
          \pf_s(t_0)^\alpha - \alpha ^\alpha \sopt(t_0)^\alpha 
                    - \alpha \pf_s(t_0)^{\alpha - 1}\pf_s(t_0)  \\ + \alpha^2 \pf_s(t_0)^{\alpha - 1}\sopt(t_0) \leq 0 \nonumber            
        \end{multline}
        Substituting $z = \frac{\pf_s(t_0)}{\sopt(t_0)}$ results in 
        \begin{align}
           \sopt(t_0)^\alpha \left( (1- \alpha) z^\alpha  + \alpha^2 z^{\alpha - 1} - \alpha ^\alpha \right) \leq 0, \nonumber
        \end{align}
        where we note that in the working case, $\sopt(t_0)>0$. Therefore, consider the polynomial
        \begin{align}
            u(z) = (1- \alpha) z^\alpha  + \alpha^2 z^{\alpha - 1} - \alpha ^\alpha. \label{eq:polynomialUofZ}
        \end{align}
        Evaluating this interval at the domain boundaries, we note that $\lim_{z\downarrow 0} u(z) = -\alpha ^\alpha$ and $\lim_{z\uparrow \infty} u(z) = -\infty$ for $\alpha > 1$. For (\ref{eq:potentialFtnIneq}) to hold, it now suffices to show that the maximum of (\ref{eq:polynomialUofZ}) does not exceed zero. To this end, we differentiate the polynomial with respect to $z$, finding
        \begin{align}
            \frac{d u}{dz}(z) = (\alpha - \alpha^2)z^{\alpha - 1} + (\alpha^3 - \alpha^2)z^{\alpha - 2}. \nonumber
        \end{align}
        Given that $z\neq 0$, this derivative assumes its unique zero in $z = \alpha$. Substituting this value into (\ref{eq:polynomialUofZ}), we find the maximum value of $u(\alpha) = 0$, thereby proving (\ref{eq:potentialFtnIneq}) for the working case.
                
        Lastly, we note that the arguments above also apply to the case where no new task arrives, but a task $j$ is completed by \opt\ or \oa\ at time $t_0$. This is so since speed profiles $\soa(t)$ and $\sopt(t)$ are unaffected by the completion of $j$, allowing us to apply the working case arguments.

This concludes the proof of Lemma~\ref{lem:potentialFunction}.
\qedsymbol

\begin{thm}\label{thm:OAcompetitiveness}
    \oa\ is $\alpha^\alpha$-competitive for DSL.    
\end{thm}    
\textit{Proof of Theorem~\ref{thm:OAcompetitiveness}: } 
We first note that if we can upper bound the competitiveness by $\alpha^\alpha$, then this bound is tight. This follows directly from Lemma~3.2 in \cite{2007BansalSpeedScalingManageEnergyAndTemperature}, where a DS instance is presented. A corresponding DSL instance given sufficiently large speed limits (e.g., speed limits $\ell_j = w_j = \left(\frac{1}{n-j}\right)^{\frac{1}{\alpha}}$), yields the same upper bound on the objective value and on the competitive ratio. For the details, we refer to~\cite{2007BansalSpeedScalingManageEnergyAndTemperature}. 

As for the upper bound on the competitive ratio, we integrate (\ref{eq:PFproperty2}) with regard to time to find that for any DSL instance $I$ and corresponding \oa\ schedule $s_{\oa}$ and optimal schedule $s_{\opt}$:
\begin{align}
    \int_t \soa(t)^\alpha + \int_t \frac{d\Phi(t)}{dt} &\leq \alpha^\alpha \int_t \sopt(t)^\alpha\label{eq:PFintegralrange} \\
    E(s_{\oa}) + \Phi(\max_j d_j) - \Phi(\min_j r_j) &\leq \alpha^\alpha E(s_{\opt}) \label{eq:PFfundamentalcalculus}\\
    E(s_{\oa}) &\leq \alpha^\alpha E(s_{\opt}) \label{eq:PFcompetitivenessBound}
\end{align}
where all integrals in \eqref{eq:PFintegralrange} are taken over the positive range $\mathbb{R}_{\geq 0}$. Furthermore, in \eqref{eq:PFfundamentalcalculus} we use the definition of the objective function, the fundamental theorem of calculus, 
and Property~\ref{item:potentialProperty1} of Lemma~\ref{lem:potentialFunction}. Finally, \eqref{eq:PFcompetitivenessBound} follows from Property~$1$ of the same lemma, i.e., from $\Phi(\max_j d_j) = \Phi(\min_j r_j) = 0$.  
\qedsymbol

\subsection{Exact Scheduling Rules}
As noted before, applying \edf\ may result in infeasible schedules for DSL instances (see Figure~\ref{fig:EDFinfeasible}). 
This section takes this observation a step further, concluding that there exists no deterministic online scheduling rule that given any speed profile corresponding to a feasible schedule can guarantee to find such a schedule. In this, we assume that the scheduling rule only becomes aware of job $j$ at its release time $r_j$, while the speed profile assumes implicit knowledge of all jobs released over the time horizon. Formally:
\begin{thm} \label{thm:noExactDeterministicOnlineSchedulingRule}
    Let $I$ be a DSL instance, and let $\pf$ be a speed profile for which there exists a feasible schedule $s$ for $I$ such that $\pf_s = \pf$. There exists no deterministic online scheduling rule that reliably finds a feasible schedule $s'$ for any such pair ($I$,$\pf$) for which $\pf_{s'} = \pf$.
\end{thm}
\begin{figure}[]
    \centering
    \begin{tabular}{c|c c c c}
         $j$ & $r_j$ & $d_j$ & $p_j$ & $\ell_j$ \\  \hline
         1 & 0 & 1 & 1 & 2 \\
         2 & 0 & 2 & 1 & 1 \\
         3 & 1 & 2 & 1 & 1 \\
         4 & 1 & 2 & 1 & 1 \\
         5 & $\frac{1}{2}$ & 1 & 1 & 2 
    \end{tabular}
    \\ \vspace{6 pt}
    
    \includegraphics{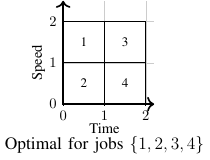}
    \includegraphics{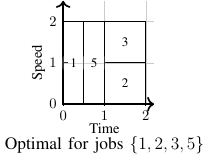}    
    \caption{Example instances for the proof that for DSL instances no online scheduling rule can reliably find a feasible schedule based on a speed profile.}
    \label{fig:NoOnlineSchedulingRuleDSL}
\end{figure}
\textit{Proof of Theorem~\ref{thm:noExactDeterministicOnlineSchedulingRule}:}    
    We prove this by providing a speed profile and two DS instances corresponding to that profile, such that for any initial scheduling decision, one of the instances cannot be feasibly scheduled. A job set with five jobs is illustrated in Figure~\ref{fig:NoOnlineSchedulingRuleDSL}. We now consider two DSL instances that are subsets of those jobs. Here, the first instance $I$ considers jobs $\mathcal{J} = \{1,2,3,4\}$ and the second instance $I'$ considers job set $\mathcal{J'} = \{1,2,3,5\}$. Note that the optimal speed profile for both instances is constant with speed 2 throughout the time horizon $[0,2]$. The respective optimal schedules are shown in the figure. Furthermore, at time $t = 0$, for both instances, only the first two jobs have been released. 

    Based on these two instances, we argue that for any scheduling decision made by an online algorithm at $t=0$, we can choose an instance such that the resulting schedule cannot feasibly follow the speed profile. First, consider the case where at $t=0$ the algorithm decides to only run Job~1. In that case, we reveal instance $I$. Limited by its speed limit, Job~2 cannot finish its processing before time $t=1$ when Jobs~3 and 4 are released.
    Next, consider the case where on time interval $[0,\frac{1}{2}]$ Job~1 and Job~2 are processed at strictly positive speed for an $\epsilon>0$ amount of time. If we now reveal instance $I'$, it is not possible anymore to follow the speed profile while completing both Jobs~1 and 5 before their respective deadlines. 
    This illustrates that no deterministic online scheduling rule can reliably follow a given speed profile, even if there does exist a feasible schedule.
\qedsymbol

In particular, the result implies that there exists no 1-consistent learning augmented online scheduling algorithm for DSL that relies on predictions of the aggregated speed profile. Even if given the optimal speed profile, the proof above indicates that there exists no deterministic scheduling rule that can reliably find an optimal schedule.

\section{Numerical experiments}
\label{sec:numericalEmpricialCompetitiveRatiosDSL}
In this section, we relate the theory developed in the previous sections to the EV scheduling application. To this end, we evaluate the discussed algorithms based on numerical experiments with real-world data.
\footnote{The code used for the simulations is available under \url{https://github.com/lwinschermann/FlowbasedOfflineChargingScheduler} (commit 7506297).}
In particular, we numerically evaluate the computational complexity of the offline algorithm in Section~\ref{sec:numerical:computational} and compare the competitive ratios of the discussed online algorithms with simulation results based on real-world data in Section~\ref{sec:numerical:competitiveness}. 

The real-world data underlying the experiments was collected at an office parking lot in Utrecht, the Netherlands between September $1^{st}$ 2022 and June $16^{th}$ 2023, resulting in a total of 13694 charging sessions\footnote{A more extensive version of the data set is available in [reference omitted for anonymization]
	.}.
A maximum of 113 charging sessions was recorded in a single day. Each recorded EV charging session is described by the EV's arrival time, departure time and the total amount of energy charged. The charging stations at the parking lot are two-plug installations that support charging with at most 11~kW or 22~kW, depending on whether one or two EVs are plugged into the same charging station. In the experiments, we chose between these two values for the individual EV-specific maximum charging rates, depending on the average power required to charge the recorded amount of energy within the EV's availability in the parking lot. Currently, there are around 250 chargers installed in the office parking lot, and this number is expected to increase to over 400 in the next few years.

\subsection{Computational efficiency of offline algorithm} \label{sec:numerical:computational}
In this section, we numerically evaluate the computational complexity of the offline algorithm \focs.
Research has shown that for some algorithms theoretical and empirical running times may differ. For example, the simplex algorithm, while theoretically exponential in running time, performs well on real-world instances. The other way around, Ellipsoid methods for linear programs are known to be polynomial on paper, while performing poorly in practice. 
This motivates us to investigate the empirical running time of the novel algorithm \focs\ based on real-world EV charging data.
We solve instances based on real-world data using a proof-of-concept implementation of \focs, and compare it to results achieved with the commercial solver Gurobi~\cite{gurobi}.

\subsubsection{Software implementation} \label{sec:numerical:software}
The proof-of-concept \focs\ implementation used for the empirical running time analysis of \focs\ has been developed specifically as a proof of concept for this paper. The code is written in python, heavily relying on the \texttt{networkx} package~\cite{2008Networkx}. This package is chosen for its user friendly flow models, and because various maximum flow algorithms are readily available within the package. An overview of the considered algorithms is provided in Table~\ref{tab:maxflowComplexity}.
\begin{table*}[]
	\caption{Complexity table for four maximum flow algorithms. Here, $N$, $M$ and $U$ are respectively the number of nodes, edges and the maximum capacity in a given flow network. The theoretical complexity is taken from \cite{2023CruzMejiaSurveyExactAlgsMaxFlowAndMinCostFlow}.}
	\centering
	\begin{tabular}{l|c|c|c}
		Maximum flow algorithm & Theoretical compl. 
		& Networkx compl. & 
		\focs\ network \\\hline 
		\texttt{shortest\_augmenting\_path} & $\mathcal{O}(N M U)$ & $\mathcal{O}(N^2 M)$ & $\mathcal{O}(n^4)$ \\ 
		\texttt{edmonds\_karp} & $\mathcal{O}(M^2N)$ or $\mathcal{O}(M^2\log U)$ & $\mathcal{O}(M^2N)$ & $\mathcal{O}(n^5)$ \\ 
		\texttt{preflow\_push} & $\mathcal{O}(MN\log (\frac{N^2}{M}))$ & $\mathcal{O}(N^2\sqrt{M})$ & $\mathcal{O}(n^3)$ \\ 
		\texttt{dinitz} & $\mathcal{O}(N^2M)$ or $\mathcal{O}(MN\log U)$ & $\mathcal{O}(N^2M)$ & $\mathcal{O}(n^4)$ \\ 
	\end{tabular} 
	\label{tab:maxflowComplexity}
\end{table*}
%
The last column of the table specifies the computational complexity of the respective algorithms when substituting the network size of the \focs\ network ($\leq 3n + 1$ nodes and $\leq 2n^2 + 2n -1$ edges, cf.~\ref{sec:flow}) into the complexity as reported by Networkx\footnote{The website was last accessed on December 19$^{th}$ 2023.} (third column in Table~\ref{tab:maxflowComplexity}). We find that \texttt{shortest\_augmenting\_path()}, \texttt{edmonds\_karp()}, \texttt{preflow\_push()} and \texttt{dinitz()} have a respective theoretical complexity of $\mathcal{O}(n^4)$, $\mathcal{O}(n^5)$, $\mathcal{O}(n^3)$ and $\mathcal{O}(n^4)$ on \focs\ networks.
Note that the implementations available in 
\texttt{networkx} likely do not strictly follow the original name-sake algorithms (\cite{1956FordFulkersonMaxFlows}, \cite{1972EdmondsKarpMaxFlowAlgs}, \cite{1974KarzanovMaxFlowsPreflows}/\cite{1988GoldbergMaxFlows} and \cite{1970DinitzMaxFlows} respectively). Since their publication, various improvements and variants have been introduced that are usually referred to by the same name.

Based on the problem constraints \eqref{eq:MIP} and the quadratic objective function \eqref{eq:energyFunctionObjective} with $\alpha = 2$, it is also possible to use the commercial solver Gurobi~\cite{gurobi} to solve DSL instances. In this paper, we use Gurobi for comparison of \focs\ with state-of-the art solvers.

\subsubsection{Experimental setup}
In this section, we describe the experimental setup used to numerically evaluate the computational efficiency of the offline algorithm \focs\ compared to the benchmark results by the commercial solver Gurobi. 

The main parameters to classify the experiments are:
\begin{itemize}
	\item instance size $n$, i.e., the number of EV charging sessions;
	\item time granularity, i.e., either 1 minute, 15 minutes or 1 hour;
	\item used maximum flow method in \focs.
\end{itemize}
For a given set of parameters for the experiment, we randomly sample $n$ sessions from the dataset described in Section~\ref{sec:numericalEmpricialCompetitiveRatiosDSL} and solve the instance using \focs\ and Gurobi. We record their CPU running times using the function \texttt{time.process\_time()} from the python package \texttt{time}. 
To get meaningful results, we repeat this process 500 times, starting at the instance sampling. The values reported in this work are the median running times over those 500 runs.
All experiments are run on an Intel Xeon E5-2630 v3 processor. 

\subsubsection{Results}
In this section, we present and analyze the results of the running time experiments. In particular, we are interested in the efficiency and therefore usability of the offline algorithm \focs\ for (large) EV parking lots. To this end, we focus on the dependency of the running time on the instance size $n$. We discuss the impact that different maximum flow algorithms have on the performance of \focs\ and compare the impact of various time granularities under invariant maximum flow methods. 

\begin{figure}
	\centering
	\includegraphics{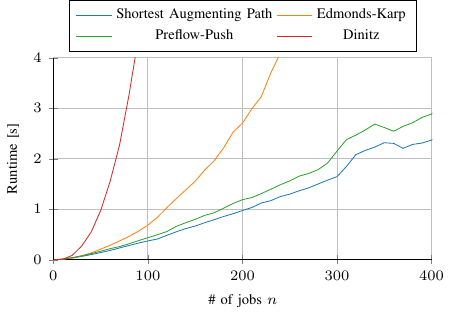}
	\caption{Running times for solving only for quarterly granularity, comparing maximum flow algorithms applied in \focs.}
	\label{fig:0123_900_full_sol}
\end{figure}

As discussed in Lemma~\ref{lem:focsfeasibleterminates}, \focs\ derives a solution in $\mathcal{O}(n^2\mu)$ time, where $\mathcal{O}(\mu)$ is the time complexity of the used maximum flow algorithm. Figure~\ref{fig:0123_900_full_sol} presents the running time of \focs\ for a full day and quarterly time granularity for four maximum flow methods. 
We see a clear difference in performance between methods that rely on augmenting paths (\texttt{shortest\_augmenting\_path()},  \texttt{edmonds\_karp()} and \texttt{dinitz()}). While the time taken to solve \focs\ using \texttt{shortest\_augmenting\_path()} seems to grow almost linearly for instance sizes up to 400, running times with \texttt{edmonds\_karp()} and \texttt{dinitz()} clearly increase non-linearly. 
Finally, the preflow-push method behaves similarly to \texttt{shortest\_augmenting\_path()}.

If we compare the above to the theoretical complexity of the various maximum flow methods embedded in \texttt{networkx} (see Table~\ref{tab:maxflowComplexity}), the most striking observation is that while \texttt{preflow\_push()} has the smallest theoretical running time on \focs\ networks ($\mathcal{O}(n^3)$), the theoretically quartic running time of \texttt{shortest\_augmenting\_path()} ($\mathcal{O}(n^4)$) outperforms the other tested methods. \texttt{Edmonds\_karp()} shows the highest theoretical complexity, and empirically is the second slowest method, overtaken only by \texttt{dinitz()}.

Overall, \texttt{shortest\_augmenting\_path()} appears to be the dominating method in our experiments with respect to running time. Therefore, for the remainder of the section we focus on results generated with \texttt{shortest\_augmenting\_path()} as subroutine in \focs. The experiments further indicate that the proper choice of maximum flow method greatly influences the usability of \focs\ for EV scheduling applications.

\begin{figure*}[t]
	\begin{center}
		\includegraphics{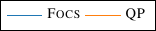}
	\end{center}
	\centering
	\subfloat[\centering 60s time granularity]{
		\includegraphics[width = 0.3\textwidth]{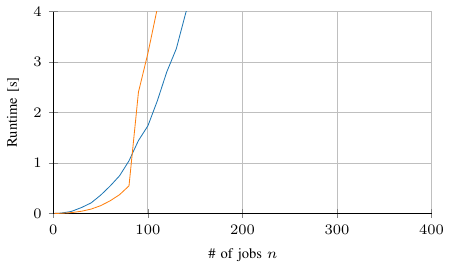}
		\label{fig:0_60_full}
	}
	\subfloat[\centering 15m time granularity]{
		\includegraphics[width=0.3 \textwidth]{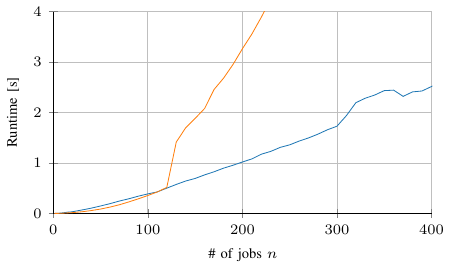}
		\label{fig:0_900_full}
	}
	\subfloat[\centering 1h time granularity]{
		\includegraphics[width = 0.3\textwidth]{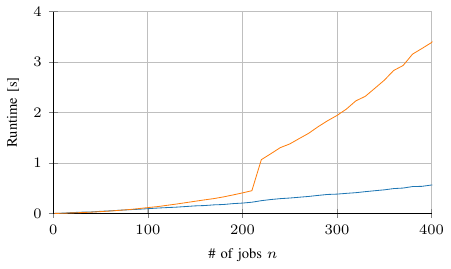}
		\label{fig:0_3600_full}
	}
	\caption{Runtimes relative to instance size, for various time granularities, solved for the whole day and using \texttt{shortest\_augmenting\_path()}.}
	\label{fig:0_full}
\end{figure*}

In the following, we compare \focs\ to the Gurobi implementation. The total running time relative to instance size $n$ for the three time granularities is depicted in Figure~\ref{fig:0_full}. 
First, we note that across subfigures, for small instances the quadratic program solved in Gurobi either outperforms \focs\ or runs at similar speed. However, for Gurobi we observe a gradually steeper increase in running time relative to instance size than for \focs. As a result, the median running time of \focs\ outperforms Gurobi starting from instance sizes of approximately 90, 120 and 80 for time granularities of 1~minute, 15~minutes and 1~hour respectively.
Furthermore, we observe a sudden steep increase in the running time of the Gurobi implementation, at instance sizes from around 90, 130 and 230 in Figures~\ref{fig:0_60_full}, \ref{fig:0_900_full} and \ref{fig:0_3600_full} respectively. This behavior possibly reflects that a certain memory threshold is reached after which solving becomes more costly for Gurobi.

Overall, as expected, a finer time granularity for a given instance size results in a larger running time for both solvers. Notably, the slope of the running time for \focs\ seems to be almost linear for the instance sizes considered, except in Figure~\ref{fig:0_60_full}. There, a non-linear increase can be observed.  
To give an indication of this growth, the median running times for \focs\ and Gurobi with instance size 400 under time granularity 1~minute are respectively 36 and 62 seconds. 

We conclude that even for a proof-of-concept implementation of \focs, in terms of running time the method is competitive with existing commercial solvers, and therefore promising for application in the field. In practice, schedules are often made in 15~minute granularity, for which \focs\ computed an optimal schedule in 2.5~seconds (median) for instance sizes up to 400. Based on the evaluation of the impact of the chosen maximum flow algorithm, we advise careful consideration of the applied maximum flow algorithm when implementing \focs\ for practical applications.

\subsection{Competitiveness of online algorithms}\label{sec:numerical:competitiveness}
As theoretical competitive ratios may be based on worst case instances that may be very unrealistic to occur in practice, we compare the competitive ratios for the online algorithms \avr\ and \oa\ presented in Section~\ref{sec:onlineDSL} with simulation results based on real-world data. We define the \emph{empirical ratios} observed in the simulations in accordance with the definition of the competitive ratio in (\ref{eq:defCompetitiveRatio}) to be $\frac{E(s^{\textsc{Alg}}(I))}{E(s^*(I))}$,
where $I$ is the considered instance, the numerator is the objective value of the considered (online) algorithm, and the denominator the objective value of the optimal solution.

For the numerical experiments, we randomly sampled 400 charging sessions and combined them into one instance. We solve this instance using \focs, \avr\ and \oa. As \focs\ is an offline solver which results in an optimal solution (see Section~\ref{sec:offlineDSL}), we use its objective value for the calculation of the empirical ratios for online algorithms \avr\ and \oa. Note that this means that \focs\ derives a schedule using perfect knowledge on all jobs, whereas \oa\ and \avr\ only become aware of jobs as they are released. We consider objective function (\ref{eq:energyFunctionObjective}) with $\alpha = 2$, and repeat the sampling and solving process 500 times, each time recording the objective value and power profile for the three algorithms. For comparison, we also evaluate a Greedy algorithm that schedules each job at its maximum speed upon release until its processing requirement is met. In practice, that corresponds to EVs arriving and charging uncontrolled at their maximum charging power until their charging demand is met, which is the default if no scheduling is applied.

\begin{figure}[t]
\centering
    \pgfplotsset{every axis legend/.append style={font={\scriptsize}}}
	\pgfplotsset{every axis legend/.append style={at={(0.5,1.03)},anchor=south}}
	\includegraphics{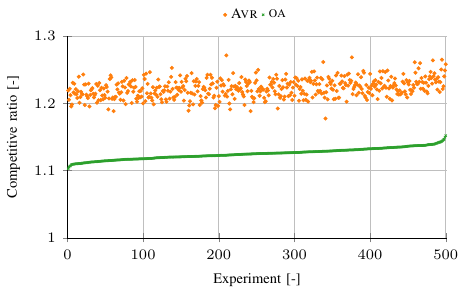}
\caption{Empirical ratios for \avr\ and \oa\ for 500 instances of size 400. Their respective theoretical bounds are 8 and 4. Results sorted in increasing order based on their empirical ratios for \oa .}
\label{fig:cr_ldc}
\end{figure}

\begin{figure}[t]
\centering
    \pgfplotsset{every axis legend/.append style={font={\scriptsize}}}
	\pgfplotsset{every axis legend/.append style={at={(0.5,1.03)},anchor=south}}
     \includegraphics[width=0.5\textwidth]{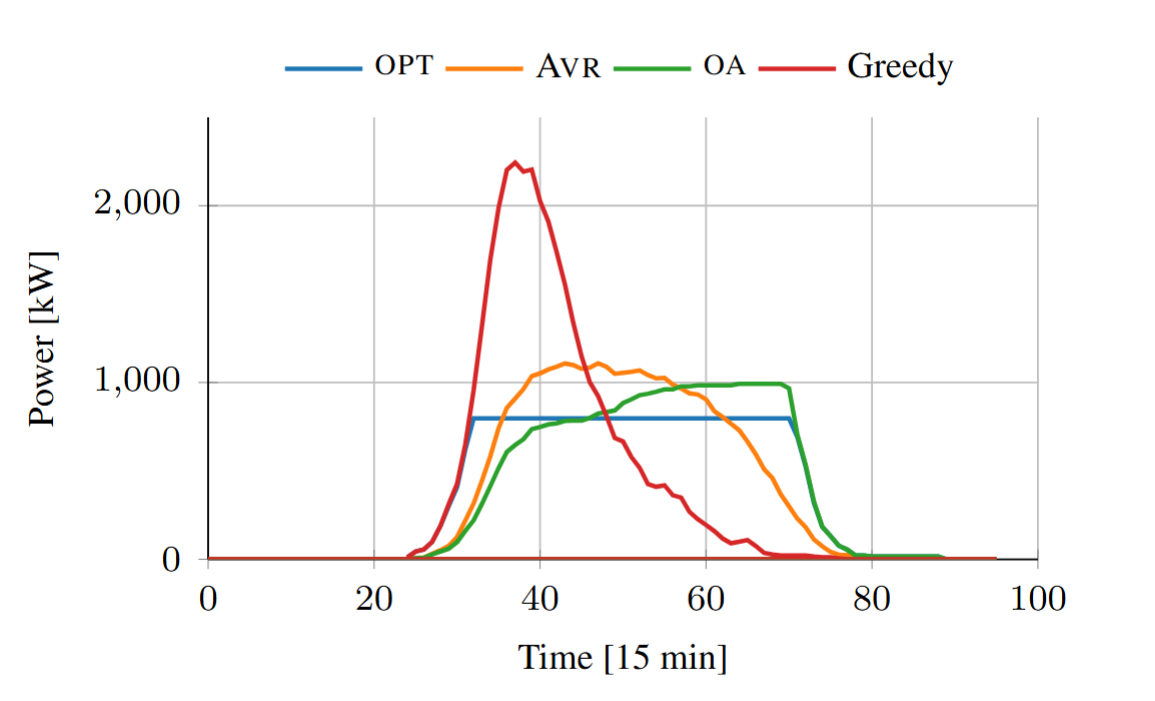}
 \caption{Aggregated speed profiles for \focs, \oa\ and \avr\ for one sampled instance based on real-world data.}
\label{fig:power_AVR_OA_OPT}
\end{figure}

The results are summarized in Figures~\ref{fig:cr_ldc} and~\ref{fig:power_AVR_OA_OPT}. 
First, Figure~\ref{fig:cr_ldc} shows the empirical ratios for \avr\ (orange) and \oa\ (green) for 500 randomly sampled instances. The instances have been sorted based on the empirical ratio for \oa, resulting in the green dots forming an increasing sequence. 
Notably, this ordering has not translated to the empirical ratios of \avr, meaning that the ordering of two DSL instances based on the objective value for \oa\ in general does not say anything about their objective values for \avr.
For \avr\ the minimum and maximum empirical ratios recorded in the experiments were 1.18 and 1.27 respectively, as opposed to the theoretical bound $2^{\alpha - 1}\alpha^{\alpha} = 8$. For \oa\ the minimum and maximum empirical ratios were 1.10 and 1.15 respectively, as opposed to the theoretically tight competitive ratio $\alpha^{\alpha} = 4$. Note that the maximum empirical ratio for \oa\ is smaller than the minimum for \avr. That implies that in terms of objective value, \oa\ dominates \avr. For comparison, the empirical ratio for Greedy assumes values between 1.68 and 2.06.

For the first out of the 500 sampled experiments, Figure~\ref{fig:power_AVR_OA_OPT} shows the four power profiles resulting from \focs\ (blue), \avr\ (orange), \oa\ (green) and Greedy (red). 
The speed profile of \avr\ is more impacted by the high simultaneity of arrivals and great variance of departure times in this particular office building \cite{2023WinschermannVoI}. 
The effect is graphically reflected by the slightly left-leaning form of the curve. Its relative smoothness can be attributed to the long dwell times and simultaneity in office parking lots. \oa, on the other hand, shifts a lot of the work to the end of the time horizon, as it is oblivious to each subsequent new arrival. From a user perspective, having a gradual charging process over the day invokes less anxiety and mistrust as opposed to charging later in the day. Furthermore, the \avr\ solution is more robust to early departure by individual EVs. This last observation falls outside of the DSL problem statement, but becomes relevant when considering e.g., charging guarantees as deterministic input to the optimization. We note that from a parking-lot perspective, both \avr\ and \oa\ clearly outperform the uncontrolled Greedy approach. From a parking lot-level perspective, a reduction of the power peak to about half the original peak is highly desirable. 

\section{Conclusion} \label{sec:conclusionDSL}
In this work, we consider an EV scheduling problem with as objective to minimize an increasing, strictly convex and differentiable function of the aggregated power profile. In particular, we relate EV scheduling to speed scaling with job-specific speed limits 
and present \focs, an exact offline algorithm that determines an optimal schedule in $\mathcal{O}(n^2 \mu)$ time where $\mathcal{O}(\mu)$ is the complexity of the used maximum flow algorithm. Numerical experiments based on real-world data show that a proof-of-concept implementation of \focs\ has a median running time of 2.5~seconds to solve instances with 400 EVs over a full day in 15 minute granularity. Given that within the energy domain it is common to consider 15~minute planning intervals and that the efficiency of the software implementation used for our experiments can still be improved, this indicates usability of \focs\ as optimization subroutine in the field. 

Next to the offline algorithm, we analyze two online algorithms and their respective competitive ratios for a class of objective functions depending on a parameter $\alpha$, where in energy applications $\alpha$ typically equals 2. Average Rate is shown to be $2^{\alpha-1}\alpha^{\alpha}$ competitive, and Optimal Available has a tight competitive ratio $\alpha^\alpha$. 
These competitive ratios match those for the classical speed scaling model that does not consider job-specific speed limits and only processes one job at a time. In our experiments with real-world EV charging data, both Average Rate and Optimal Available approximate the optimal offline solution with a factor of less than 1.3. This performance is significantly better than the uncontrolled default often used in practice where EVs charge at their maximum power starting from their arrival until their energy charging demand is met. 

Future work may investigate additional problem constraints such as global power limits. Furthermore, numerical experiments are of interest, especially their integration with control strategies such as model predictive control or fill-level algorithms. Lastly, given that optimal schedules are not necessarily unique, scheduling rules resulting in a robust output should be explored.

\bibliographystyle{IEEEtran}
\bibliography{literature_references_and_summaries}

\end{document}